\documentclass{compositio}
\usepackage{amssymb}
\usepackage{amsfonts}
\usepackage{amstext}
\usepackage{amsmath}
\usepackage{amsthm}
\input xypic
\xyoption{all}

\newtheorem{thm}{Theorem}[section]
\newtheorem{defn}[thm]{Definition}
\newtheorem{prop}[thm]{Proposition}
\newtheorem{lem}[thm]{Lemma}
\newtheorem{cor}[thm]{Corollary}

\newtheorem{proposition}[thm]{Proposition}
\newtheorem{lemma}[thm]{Lemma}

\newtheorem{remark}[thm]{Remark}
\newtheorem{fact}[thm]{Fact}

\def\AA{\mathbb A}

\def\CC{\mathbb C}

\def\NN{\mathbb N}

\def\PP{\mathbb P}
\def\QQ{\mathbb Q}

\def\ZZ{\mathbb Z}

\def\aa{\mathcal A}

\def\cc{\mathcal C}
\def\dd{\mathcal D}

\def\ff{\mathcal F}
\def\gg{\mathcal G}

\def\ll{\mathcal L}
\def\mm{\mathcal M}

\def\pp{\mathcal P}

\def\Hom{{\mathrm{Hom}}}
\def\Id{{\mathrm{Id}}}
\def\CLD{{\mathrm{C}}_*}
\def\indehm{{\mathrm{Ind} \mathrm{EHM}}}

\def\End{{\mathrm{End}}}
\def\sing{{\mathrm{Sing}}}
\def\Aff{{\mathrm{Aff}}}

\def\spec{{\mathrm{Spec}}}
\def\pf{{\noindent \textit{Proof. }}}
\def\endpf{$\square$ \bigskip}

\def\forg{f \hspace{-2pt} f}
\def\BAR{\mathrm{Bar}}
\def\prim{\mathrm{Prim}}
\def\fin{\mathrm{Fin}}
\def\opp{\mathrm{opp}}
\def\pairs{\mathrm{Pairs}}
\def\tensor{\otimes }
\def\Zmod{\ZZ\mathrm{-mod}}
\def\HomS{\Hom_\dd(\bullet, S)}

\def\Ch{\mathrm{Ch}}
\def\Ab{\mathcal{A}b}
\def\im{\mathrm{im}}
\def\BART{\widetilde{ \mathrm{Bar}}}
\def\sets{\mathrm{Sets}}

\def\Ab{\mathrm{Ab}}
\def\id{\mathrm{Id}}
\def\ehm{\mathrm{EHM}}
\def\dim{\mathrm{dim}}

\begin{document}

\title[A Construction of a D.G. Lie Algebra in E.H.M.]{A Construction of a Differential Graded Lie Algebra in the Category of Effective Homological Motives}
\author{Kaj M. Gartz}
\email{gartz@math.uchicago.edu}
\address{Mathematics Department\\University of Chicago\\5734 S. University Ave, Chicago IL 60637\\USA }

\classification{14F35, 14F42, 57T99}
\date{}

\begin{abstract}

This text gives a construction of a differential graded Lie algebra
in Nori's category of effective homological motives. In fact the
construction works in more a general setting than that of an Abelian
category. This allows us to give the rational homotopy Lie algebra of
a $1$-connected space a motivic structure. As a consequence the
rational homotopy Lie algebra inherits a mixed Hodge structure and
Galois module structure.

\end{abstract}

\maketitle

\section{Introduction}

The goal of this paper is to show that the rational homotopy Lie
algebra of a variety (defined in Section \ref{sec:RHBD})  carries a
motivic structure in the sense of Nori. This is formalized as
follows. If $k \hookrightarrow \CC$ is an embedding of a field $k$
into the complex numbers then Nori (see \cite{br} and \cite{le})
defines an Abelian tensor category of effective homological motives
over $k$, called $\indehm$. Let $\mm = \Ch(\indehm)$ the Abelian
category of chain complexes of motives. To say that $X$ carries a
motivic structure means that $X$ is an object in $\mm$ and that
$H_*(X)$ is an object of $\ehm$.

Nori uses his Basic Lemma \cite{nori} to construct a functor which
assigns to any variety a complex of motives. A simpler case of this
construction occurs when one considers the category of affine
varieties.  We denote this functor by $X \mapsto \CLD(X)$.

The restriction to affine varieties is not too severe thanks to
``Jouanolou's Trick'' (see \cite{jo}): For any quasi-projective
variety $X$ over a field $k$, there exists an affine variety $X'$
over $k$ and a morphism $X' \to X$ which is  a Zariski locally
trivial fibration with fibers isomorphic to $\AA^n$. Since the affine
fibers are contractible, $X'$ is homotopy equivalent to $X$.

Nori also provides realization functors from $\indehm$ to the
category of mixed Hodge structures and to the category of Galois
representations. Let $\mm$ be the category of chain complexes in
$\indehm$. There is a forgetful functor, $\forg$ from $\indehm$ to
Abelian groups (and thus a forgetful functor from $\mm$ to chain
complexes of Abelian groups), such that
    \begin{equation*}
        H_*(\forg \CLD(X)) \cong H_*(X(\CC), \ZZ),
    \end{equation*}
where $H_*(X(\CC), \ZZ)$ denotes the singular homology of the
topological space $X(\CC)$.

In this paper we provide a functor $\pp_F$ which associates to any
affine variety, $X$, a differential graded Lie algebra (d.g.l.),
$\pp_F(X)$ in $\mm$, such that if the variety is simply connected
then the homology of the d.g.l. computes the rational homotopy Lie
algebra, that is
    \begin{equation*}
    H_*( \pp_F) \cong \pi_{*+1}(X)\otimes \QQ
    \end{equation*}
as graded Lie algebras. The known topological techniques for
producing such a Lie algebra all use  a coalgebra analogous to that
of singular chains with the Alexander-Whitney map or the wedge
product of forms.

This functor, $\pp_F$, comes from an underlying combinatorial
construction on geometric objects, and is perhaps the simplest Lie
algebra one can construct from a simplicial object in an Abelian
category which is representable (i.e. it can be extended to the
category of finite sets with all morphisms, not just order preserving
maps.) Applying the construction with singular chains to any
topological space yields a new method of producing a d.g.l. which
computes the rational homotopy Lie algebra.

This functorial construction has several applications. First, since
Nori's category of motives has realization functors, our functor
gives rise to a Galois structure on the homotopy Lie algebra after
tensoring with $\QQ_p$. Previously there was the mixed Hodge
structure provided by Hain, but since he uses the graded commutative
algebra structure of the deRham complex, his techniques were not
easily adapted to the case of \'etale cohomology and Galois
representations. Second, we pave the way for a motivic study of
extensions arising from higher rational homotopy theory -- analogous
to Hain's work with extensions of mixed Hodge structures arising from
the fundamental group \cite{hain2}.

The serious difficulty in showing that the rational homotopy Lie
algebra of a variety has a motivic structure reduces to the fact that
although there is a motivic analog of the Eilenberg-Zilber map
    \begin{equation*}
        \CLD(X) \otimes \CLD(Y) \to \CLD( X \times Y),
    \end{equation*}
which is a quasi-isomorphism, there is no map in the other direction
corresponding to the Alexander-Whitney map. Without the
Alexander-Whitney map there is no coalgebra structure on motivic
chains.
%
It is unreasonable to expect a splitting of this map (like the
Alexander-Whitney map) in a category of motives (otherwise the
splitting forces non-zero extension classes to vanish). Thus
$\CLD(X)$ does not have the structure of a coalgebra, in contrast to
the topological case, where one has maps on singular chains
    \begin{equation*}
    \sing(X) \xrightarrow{\sing(\Delta)} \sing(X \times X)
    \xrightarrow{\mathrm{A-W}} \sing(X) \otimes \sing(X).
    \end{equation*}
    Put another way -- it is clear how
to take two cycles and produce a cycle on the product, but not vice
versa. Thus much of the machinery of algebraic topology (in our case
in particular anything related to differential graded Hopf algebras)
can not be reproduced in the motivic category without much work.

Although  we have no coalgebra, we still have a map $
\CLD(X)^{\tensor n} \to \CLD(X^n) $ which is $\Sigma_n$ equivariant.
In some ways the key to this paper is that this
$\Sigma_n$--equivariance is enough to get a Lie algebra and that we
don't need (or have) a coalgebra structure.


One crucial point is that when $V$ is a graded vector space,
considering the tensor algebra $T(V)$ as a differential graded Hopf
algebra, the primitive Lie algebra in $T(V)$ coincides with the free
Lie algebra on $V$, $L(V)$ -- and this can be obtained as the {\it
image} of a projector defined in terms of the $\Sigma_n$ action on
$\CLD(X^n)$. This is crucial because defining this Lie algebra as
primitives of a co-multiplication map would require a map $\CLD(X)
\to \CLD(X) \tensor \CLD(X)$ which does not exist since there is no
map $\CLD( X \times X) \to \CLD(X) \tensor \CLD(X)$. The construction
of $\pp_F$ is very simple and combinatorial in nature. It has not
been previously discussed before probably because: 1. Algebraic
topologists would have no need to make this construction with the
presence of the Alexander Whitney map (and thus the rich structure of
Hopf algebras \cite{milnor}) 2. This construction lives outside of
the simplicial category as we need to consider {\it all} maps between
finite sets 3. Algebraic geometers such as \cite{hain1} would make
use of commutative algebras such as those arising from differential
forms.

Several others have worked on related problems: Cushman's thesis,
which produced a motivic structure on the group ring of the
fundamental group completed at the augmentation ideal \cite{cush},
provided a crucial idea that the geometric cobar complex could be
used in Nori's category to compute the homology of the loop space. In
the setting of mixed Hodge structures, building on the work of Chen,
Hain showed that for any variety over $\CC$, the Malcev completion of
the fundamental group and the rational homotopy Lie algebra (in the
simply connected case) carry mixed Hodge structures
\cite{hain2,hain1}. This work was further pursued  by Wojtkowiak, who
stressed the usefulness of using cosimplicial varieties \cite{wo}.

The paper is organized as follows. In section 2, we lay out the
categorical framework which we find most convenient to describe the
d.g.l. and develop the algebraic identities in the category of
Abelian groups which are used to deduce identities in a more general
additive category $\aa$. In section 3 we give the construction of the
d.g.l. in $\aa$. In section 4 we specialize to the category of
effective homological motives. We summarize without proof many of the
properties of Nori's construction of motives. This section reduces
the problem to the case of singular chains, and thus to a problem in
topology. In section 5 we compare this construction with the work of
Hain to show that the forgetful functor applied to the d.g.l. in
$\mm$ gives a complex which computes the rational homotopy of a
$1$-connected variety.

The key mathematics in this work is in Section 3, in particular
Section 3.2, where we show $\pp_F$ is closed under its differential
$F(f_n)$. Note that in diagram \eqref{eq:digweed} $\phi_n$ is not
$f_n$! Much of the insight needed to prove this theorem was obtained
by constructing this $\phi_n$ ``by hand". This was related to
studying base-point dependency in the study of extensions of motives
arising from the fundamental group of $\PP^1- \{ 0, 1, \infty \}$
\cite{deligne}.

The author would like to thank Madhav Nori for his help, inspiration,
and support.

\section{Preliminaries}

\subsection{Conventions and Definitions}

Let $\fin$ be the category whose objects are finite sets and whose
morphisms are all set maps (not just order preserving maps). Very
contrary to the simplicial literature, we let $[n]$ denote a finite
set with $n$ elements labelled $1, 2, \ldots , n$.

The symmetric group on $n$ letters, $\Sigma_n$, is the set of all
bijections from $[n]$ to $[n]$. Since elements of $\Sigma_n$ are
maps, we think of $\Sigma_n$ as acting on the left on $[n]$. Given a
functor $G: \fin^\opp \to \cc$ to some category $\cc$, the left
$\Sigma_n$ action on $[n]$ becomes a right $\Sigma_n$ action on
$G([n])$.

Often we will want to take linear combinations of maps of objects
that have no Abelian group structure. To do this, given a category
$\cc$ we construct a new category $\dd_\cc$ (or just $\dd$ if $\cc$
is understood). The objects of $\dd$ are the same as those of $\cc$
but we define
        \[ \Hom_\dd(X,Y) = \ZZ \Hom_\cc(X,Y). \]
This gives an additive category where one can make sense of the
notion of the composite of two maps being zero.  We may also tensor
with $\QQ$ to give $\Hom_\dd$ the structure of a vector space. By
definition an additive category has a zero object and a product
structure.

A functor $G: \cc \to \cc'$ induces a functor (abusively called $G$)
$G: \dd_\cc \to \dd_{\cc'}$. Furthermore, if $G: \cc \to \aa$ is a
functor to an additive category, then the functor automatically
extends uniquely to a functor $G: \dd_\cc \to \aa$. We will be
concerned with several functors in this paper:

\subsection{Four Functors}

We will be concerned with four functors on the category $\fin$.

\begin{enumerate}
\item  The functor represented by a set $S$.

Let $S$ be any set. $S$ determines a functor (the ``functor of
$[n]$-valued points of $S$'')
    \begin{equation}
            S:\fin^\opp \to \sets
    \end{equation}
    given by $S(\bullet) = \Hom_\fin(\bullet, S)$.
Notice that $S([n])=\Hom_{\mathrm{Sets}}([n],S)$ may be identified
with the product $S^n$ ($n$-tuples of elements in $S$). A morphism of
finite sets $f: [m] \to [n]$ induces a map $S^n \to S^m$ by pullback
of coordinates.

Since $S(\bullet)$ is a  functor, $S([n])=S^n$ has a right action of
$\Sigma_n$.
For example if $n \geq 3$, $\sigma_{(123)} \in \Sigma_n$ acts
on $S^n$:
    \begin{equation*}
        (x_1, x_{2}, \ldots, x_{n-1}, x_n) \sigma =
         (x_2, x_3, x_1 , \ldots , x_n)
    \end{equation*}
We will apply this functor when $S$ itself is a finite set.

\item The functor represented by a graded $\ZZ$-module on a set.

A generalization of the functor above is to take the functor
$\Hom_{\dd_\fin}(\bullet,S)$,
    \begin{equation}
            \Hom_{\dd_\fin}(\bullet,S):  \fin^\opp \to \text{Abelian Groups}
    \end{equation}
Then $\Hom_{\dd_\fin}([n],S)$ may be identified with $\ZZ[S]^{\tensor
n}$, the $n$-th tensor power of the free module on $S$. If the
elements $S$ were to be given a grading then $\ZZ[S]^{\tensor n}$
would be a graded $\ZZ$ module. We could also tensor with a field $k$
to obtain a vector space $k[S]^{\tensor n}$.

$\Sigma_n$ acts on the $n$-th tensor power of $\ZZ[S]$ via the action
on $S$. For example notice that the map $[2] \to [1]$ induces the
linear map which on the basis $S$, sends $s_i \mapsto s_i \tensor
s_i$.

The fact that $V$ may be a graded vector space forces the symmetric
group action to acquire a sign, e.g. $\sigma_{(12)}$ acts on $v
\otimes w$ by
    \begin{equation*}
        \sigma( v \otimes w) = (-1)^{|v| \, |w|} w \otimes v.
    \end{equation*}
The graded rule from homological algebra is that anytime you pass one
symbol $X$ by another $Y$, you should multiply by $(-1)^{|X|\,|Y|}$,
and this convention is a reflection of that one.

\item The functor represented by a pair.

This is a generalization of the first functor in that first we
generalize from sets to topological spaces and then from topological
spaces to pairs. For a space $X$, throughout the rest of the paper we
fix a base point $p\in X$. We will write $X$ for the pair $(X,p)$,
unless otherwise indicated. Again we have
    \begin{equation}
        X^\bullet : \fin^\opp \to \pairs
    \end{equation}
    given by $[n] \mapsto X^n$
where $X^n$ is the topologist's product  $(X^n, p \times X^{n-1} \cup
\ldots \cup X^{n-1} \times p )$ (which is homotopy equivalent to the
smash $X^{\wedge n}$). Again, the morphisms are given by pull-back of
coordinates. For example the map $[2] \to [1]$ induces the diagonal
map $\Delta: X \to X \times X$.

\item The functor $F: \fin^\opp \to \aa$.

Let our initial hypothesis be that $\aa$ is an additive category with
tensor product $ \tensor : \aa \times \aa \to \aa$ that satisfies the
usual associativity and commutativity constraints (see \cite{dm}).
Also we assume there exists an identity object $1_\aa$ for the tensor
product. In the application to motives we will take $\aa = \mm$.


We also now assume the existence of a functor $F: \fin^\opp \to \aa$.
Since $\aa$ is additive, $F$ extends to a functor from
$\dd_\fin^\opp$ to $\aa$.
In application $F$ will be $F_X$ the functor determined by letting
$X$ be the pair of an affine variety and base point $p$ and taking
$F_X= \CLD \circ X^\bullet: \fin^\opp \to \Ch(\mm)$.

\end{enumerate}

\subsection{The Disjoint Union and Product Structure} \label{sec:axioms}

Given two finite sets $[p]$ and $[q]$ we can form their disjoint
union to get a set with $p+q$ elements. That is to say, we have a
functor
    \begin{equation}
        \amalg: \fin \times \fin \to \fin
    \end{equation}
If $f \times f': S \times S' \to T \times T'$ is a morphism in $\fin
\times \fin$, then we write $f\amalg f': S \amalg S' \to T \amalg T'$
for the functor $\amalg$ applied to $f \times f'$.

Let us consider the four functors mentioned above and their
compatibility with $\amalg$.

\begin{enumerate}
\item Sets. For a fixed set $S$, the following diagram obviously
commutes.
    \begin{equation}
        \xymatrix{
        \fin \times \fin \ar@<-1ex>[d]_{S \times S} \ar[r]^\amalg & \fin \ar[d]^S \\
        \sets \times \sets \ar[r]^\prod  & \sets
        }
    \end{equation}
That is, $\prod \circ (S \times S) = S \circ \amalg$.

\item Free $\ZZ$ modules with a chosen basis $S$. The following diagram
commutes:
    \begin{equation} \label{eq:vectortensoramalg}
        \xymatrix{
        \fin \times \fin \ar@<-1ex>[d]_{\HomS \times \HomS}  \ar[r]^\amalg & \fin \ar[d]^\HomS \\
        \Zmod \times \Zmod \ar[r]^{\phantom{aaa}\tensor} & \Zmod
        }
    \end{equation}
So we have $\tensor \circ (\HomS \times \HomS) = \HomS \circ \amalg$.

\item For the functor $X^\bullet$ we have the cartesian product
completing the square:
    \begin{equation}
        \xymatrix{
        \fin \times \fin \ar@<-1ex>[d]_{X^\bullet \times X^\bullet} \ar[r]^\amalg & \fin \ar[d]^X \\
        \pairs \times \pairs \ar[r]^\times & \pairs
        }
    \end{equation}
and similarly $\times \circ (X^\bullet \times X^\bullet) = X^\bullet
\circ \amalg$.

\item For the functor $F: \fin^\opp \to \aa$, we do {\bf not} require
that
    \begin{equation}
        \xymatrix{
        \fin \times \fin \ar[d]_{F \times F} \ar[r]^\amalg & \fin \ar[d]^F \\
        \aa \times \aa \ar[r]^\tensor & \aa
        }
    \end{equation}
commutes.

We weaken this and only  require that we have a natural
transformation of functors
    \begin{equation}
        N: \tensor \circ (F \times F) \to F \circ \amalg.
    \end{equation}
In practice this will be Nori's motivic version of the
Eilenberg--Zilber map. \end{enumerate}

Let us further assume that $F$ satisfies the following axioms.
Together with the axioms in Section \ref{sec:propsofF}, these are the
axioms which Nori advocates as a replacement for $E_\infty$
coalgebras. We have not explored the relationship with the work of
Segal  (\cite{se}) who proposed similar axioms in the topological
setting.

For the natural transformation on pairs, for all $S, T \in \fin$,
    \begin{equation}\begin{split}
        \otimes \circ (F \times F) \xrightarrow{N} F\circ \amalg \\
        FS \otimes FT \xrightarrow{N(S,T)} F( S \amalg T)
        \end{split}
    \end{equation}
 satisfies
    \begin{enumerate}
        \item If $S=\emptyset$ then $F(S)= 1_\aa$.
        \item For all $S \in \fin$ we have both
                \begin{equation}\begin{split}
                    F \emptyset \otimes FS \xrightarrow{N(\emptyset,S)} FS \\
                    F S \otimes F \emptyset \xrightarrow{N(S,\emptyset)} FS \end{split}
                \end{equation}
                are canonical identifications in $\aa$
        \item (commutativity) The switch map $\sigma: T \amalg S \to S \amalg T$
        induces a map \[F(\sigma):  F(S \amalg T) \to F(T \amalg
        S)\] making the following diagram commute:
            \begin{equation}\label{commutativity}
            \xymatrix{
                FS \otimes FT \ar[rr]^{N(S,T)} \ar[d] & & F(S \amalg T)
                \ar[d]^{F(\sigma)} \\
                FT \otimes FS \ar[rr]^{N(T,S)} & & F( T \amalg S)
            }
            \end{equation}
        \item (associativity) Let $R, S, T \in \fin$ then the following diagram
        commutes
            \begin{equation}\label{associativity}
            \xymatrix{
            (FR \otimes FS) \otimes FT \ar[rr]^{N(R,S) \otimes {1}} \ar[dd]& & F(R \amalg S) \otimes FT \ar[dr]^{N(R\amalg S, T)} & \\
                    &   &    & F(R \amalg S \amalg T) \\
            FR \otimes (FS \otimes FT) \ar[rr]^{{1}\otimes N(S,T)}& & FR \otimes F(S \amalg T) \ar[ur]^{N(R, S \amalg T)} &
            }
            \end{equation}
    \end{enumerate}

\begin{prop}[Corollary to the Axioms] \label{cortoaxioms}
    \begin{enumerate}
        \item  $F([p]) \tensor F([q]) \xrightarrow{N([p],[q])}
        F([p+q])$ is $\Sigma_p\times \Sigma_q$-equivariant.
        \item $(N([n-1],[1])) \circ
        (N([n-2],[1])\tensor {[1]}) \circ  \cdots  \circ ( N([1],[1]) \tensor
        {1} \tensor \cdots \tensor  {1}) :F([1])^{\tensor n}
        \rightarrow F([n])$ is $\Sigma_n$ -- equivariant.
    \end{enumerate}
\end{prop}
\pf    Both i and ii are consequences of $N$ being a natural
transformation.

\endpf

%
\subsection{Useful Algebraic Identities}\label{sec:algident}

The goal of this section is to observe some identities in the group
ring $\ZZ[\Sigma_n]$. The identities will be applied by considering
    \begin{equation}
        \ZZ[\Sigma_n] \hookrightarrow \Hom_{\dd_\fin}([n],[n]).
    \end{equation}
All vector spaces are taken over a field of characteristic zero.

We begin by describing some properties of Lie algebras in terms of
elements in $\ZZ[\Sigma_n]$.

For each $n \in \NN$ define
    \begin{equation}
        s_n=  (1- \sigma_{(12)}) \cdots (1- \sigma_{(12 \ldots n-1)})
        (1- \sigma_{(12 \ldots n)})\in \ZZ[\Sigma_n].
    \end{equation}
    If $n=1$ let $s_n= 1$.
If $V$ is a vector space (concentrated in degree $0$), then right
multiplication by $s_n$ above, $R_{s_n}$, acts on $V^{\otimes n}$
exactly as the Lie bracket:
    \begin{equation}
        \xymatrix{
            V^{\otimes n} \ar[r]^{R_{s_n}} & V^{\otimes n} \\
            v_1 \otimes \ldots \otimes v_n \ar@{|->}[r] &
            [\ldots[v_1,v_2] \ldots , v_n]
        }
    \end{equation}

For each $n\in \NN$ define $w_n \in \ZZ[\Sigma_n]$ by
    \begin{equation}\label{wndef}
        w_n= (1 + \sigma_{(12)})\ldots (1 + (-1)^{n-1}\sigma_{(123 \ldots n-1)}) \ldots
        (1 + (-1)^{n} \sigma_{(123 \ldots n)})
    \end{equation}
and if $n=1$, let $w_n= 1$.
$w_n$ is $s_n$ twisted by the sign representation, $\epsilon:
\Sigma_n \to \{ \pm 1 \}$:
    \begin{equation}
        w_n=\prod_{i=2}^{n}(1- \epsilon(\sigma_{(12 \ldots i)})\sigma_{(12 \ldots
        i)}) \in \ZZ[\Sigma_n].
    \end{equation}

If $V$ is a vector space concentrated in odd degree, then $w_n$ acts
on the right by the Lie bracket:
    \begin{equation}\label{wnisbracket}
            R_{w_n}( v_1 \otimes \ldots \otimes v_n) = ( v_1 \otimes \ldots \otimes v_n)w_n =
            [\ldots  [v_1, v_2]  \ldots ,v_n]
    \end{equation}

Let $\ll(V)$ denote the free Lie algebra on $V$. If $V$ is graded
then $\ll(V)$ is a graded Lie algebra. We write $\ll(V) =
\bigoplus_{i=1}^\infty L_i(V)$ where $L_i$ is the vector space of
$i$-fold iterated brackets of elements of $V$.

Henceforth let  $V$ be concentrated in degree one, then $L_i(V) =
V^{\otimes i}w_i$, so that $\ll(V) = \bigoplus_{i=1}^{\infty}
V^{\otimes i}w_i$. Let $T(V)$ be the tensor algebra on $V$. $T(V)$
can be identified with the universal enveloping algebra of the free
Lie algebra on $V$, giving $T(V)$ the structure of a Hopf algebra
(see \cite{milnor} and  \cite{q}). Then using the coalgebra structure
on $T(V)$, we have that $\ll(V)$ is the primitive Lie algebra of
primitive elements of $T(V)$. Now the $w_i$ assemble to give a linear
map, $w : T(V) \to T(V)$, where $w$ acts on the degree $n$ terms of
$T(V)$ as right multiplication by $w_n$. Thus the image of $w$ is
$\ll(V)$. Where normally $\ll(V)$ is thought of as the kernel of the
reduced co-multiplication induced from $U(L) \to U(L) \tensor U(L)$
(i.e. the Lie algebra of primitives) we have identified $\ll(V)$ with
the image of $w : T(V) \to T(V)$.

Notice that :
    \begin{equation}
        s_n^2=n s_n.
    \end{equation}

The analog for $w_n$ is the following Proposition.
    \begin{proposition} \label{wn2wn}
    \begin{equation}\label{wordformula}
        w_n^2=n w_n
    \end{equation}
    \end{proposition}
    \pf
        The proof of the following Lemma is used in Quillen's proof of the
        Poincar\'e-Birkhoff-Witt theorem and is contained in his paper \cite{q}:

        \begin{lem}\label{lemmaquillen}
            The map $\rho: T(V) \to \ll(V)$ given by
            \begin{equation}
                \rho( x_n \otimes x_{n-1} \otimes \ldots \otimes x_1) = \begin{cases}
                {\frac{1}{n}} \mathrm{ad}_{x_n} \mathrm{ad}_{x_{n-1}} \ldots
                \mathrm{ad}_{x_2}(x_1)& \text{ if $n>0$,}  \\ 0& \text{ if $n=0$} \end{cases}
            \end{equation}
            is a left inverse for the map $\ll(V) \to T(V)$
        \end{lem}

  $\rho$ is almost the linear map given by right multiplication by $\frac{1}{n}
  w_n$, since $\rho$ sends $ v_1 \otimes \ldots \otimes v_n$ to
  $\frac{1}{n} [x_1, [x_2, \ldots [x_{n-1},x_n]\ldots ]$ while
  $\frac{1}{n}R_{w_n}$ sends $v_1 \otimes \ldots \otimes v_n$ to
   \[ [ \ldots [x_1,x_2], x_3], \ldots, x_n]. \] Clearly
  $\frac{1}{n}R_{w_n}$ is  also a left inverse for the map $\ll(V)
  \to T(V)$.

     Given $l\in
    \ll(V)$, since $\ll(V)$ is the free Lie algebra on a vector
    space $V$,
    $l$ is a linear combination (after using the Jacobi identity) of terms of the form
    \[ [ \ldots [x_1,x_2], x_3], \ldots x_n]. \] Considering this as an element of
    $T(V)$, $l$ can be expressed as a linear combination of terms like $R_{w_n}(x_1\otimes \ldots \otimes x_n)$.
    The Lemma states that applying ${\frac{1}{n}}R_{w_n}$ to these
    elements is the identity:
        \begin{equation}
            (x_1 \otimes \ldots \otimes x_n)  w_n \circ {\frac{1}{n}} w_n  =  (x_1 \otimes
            \ldots \otimes x_n)w_n.
        \end{equation}
    Multiplying both sides by $n$ gives $w_n^2=n w_n$. \endpf

Notice that when we identify $\ll(V)$ with $\text{im}(w) \subset
T(V)$, we have that when restricting $w$ to degree $n$, $w^2= n w$.
Thus the following sequence is exact
    \begin{equation} \xymatrix{
        0 \ar[r] & (R_{w_n} V^{\otimes n}) \ar@{^{(}->}[r] &  {V^{\otimes
        n}}
        \ar@{^{(}->}[rr]^{R_{n - w_n}}& & V^{\otimes n}. }
    \end{equation}
This may be rewritten; in degree $n$ we have:
    \begin{equation}\label{exact}
        \xymatrix{
        0 \ar[r] &  L_n(V) \ar@{^{(}->}[r]^i & (TV)_n
        \ar[rr]^{R_{n - w_n}} & & (TV)_n
        }
    \end{equation}
Thus $R_{w_n}: (TV)_n \twoheadrightarrow L_n(V)$ is a splitting for
$L_n(V) \hookrightarrow (TV)_n$.

Let us recall some definitions. Let $R$ be a graded algebra over a
field. A graded algebra derivation of degree $|d|$ is a linear map
$d: R \to R$ which satisfies
    \begin{equation}\label{derivation}
        d(u \cdot v)=du \cdot v + (-1)^{|u| \, |d|} u \cdot dv
    \end{equation}
This implies that when we consider $R$, with the bracket given by
$[u,v] = uv - (-1)^{|u||v|} vu$,  $d$ is a Lie derivation:
    \begin{equation}\label{liealgderiv}
        d[u,v]= [du,v] + (-1)^{|u||d|}[u, dv].
    \end{equation}

Now fix a basis $v_1, \ldots v_m$ for a vector space $V$ concentrated
in degree 1. We will consider the set $S= \{ v_1, \ldots, v_m \}$ and
apply the functor $\Hom_{\dd_\fin}([n],S)\tensor k$ in the next
section. Let $Dx =n x$ if $x \in L_n(V)$. Then a straightforward
computation shows that $D: \ll(V) \to \ll(V)$ is a derivation of
degree 0.

For $w: TV \to TV$ defined above, we see that $D = w$ when restricted
to  $\ll(V)$. That is $w$ is {\it not} an algebra derivation of $TV$,
though it is the Lie algebra derivation given by $D$ on $\ll(V)
\subset TV$.

Before defining our second useful derivation we need the following
    \begin{lemma} \label{extendderiv}
        \begin{enumerate}
        \item[a.] Let $V$ be a graded vector space concentrated in degree $1$. Suppose that
        $\alpha: V \to TV$ is a linear map ``of pure degree $k$,'' (i.e. so that $V \to V^{\otimes k+1}$). Then  $\alpha$ extends
        to a graded derivation $d: TV \to TV$.
        \item[b.] If $\alpha(V) \subset \ll(V)$, then $d$ restricts to
        a Lie derivation, and thus $d( \ll(V)) \subset \ll(V)$.
        \end{enumerate}
    \end{lemma}
    \pf
        a. To give a derivation on an algebra we need only define
        $d$ on generators. Since $TV$ is the free algebra on $V$,
        $d$ extends to products by
            \begin{equation}
                d(v_1 \otimes \ldots \otimes v_n) =
                \sum_{i=1}^n(-1)^{(i-1)k} v_1 \otimes \ldots \otimes
                v_{i-1}\otimes \alpha(v_i) \otimes v_{i+1}
                \otimes \ldots \otimes v_n.
            \end{equation}

        b. An algebra derivation gives a Lie algebra derivation as
        in \eqref{liealgderiv}. Now $d$ takes the generators of
        $\ll(V)$ to $\ll(V)$ by hypothesis. To show that
        $d(\ll(V)) \subset \ll(V)$ proceed by induction on bracket
        length: if
            \begin{equation}
                d ( \oplus_{i=1}^{k} L_i(V)) \subset \ll(V)
            \end{equation}
            then write (using Jacobi) any element $x\in L_{k+1}$ as a linear
            combination of elements of the form $[y,z]$ where $y$ and
            $z$ are of bracket length at most $k$. Then applying
            formula \eqref{liealgderiv} we see $d[y,z]$ is written
            as a sum of $[dy,z]$ and $[y,dz]$ which are both in
            $\ll(V)$ by induction.
        \endpf

We now introduce the second derivation (still using the basis $v_1,
\ldots, v_m$ concentrated in degree one). Define
    \begin{equation} \label{geoderiv}
        \alpha: V \to V \otimes V \subset \ll(V)\\
    \end{equation}
by  $v_i \mapsto v_i \otimes v_i$ and extend by linearity.  By Lemma
\ref{extendderiv}, $\alpha$ extends to
    \begin{equation}
        d: \ll(V) \to \ll(V).
    \end{equation}
Though the constructions have so far taken place over a field of
characteristic zero, one might wish to make similar constructions
over $\ZZ$. That is, we could consider algebras and Lie algebras over
$\ZZ$ rather than over a field. Notice that if we were to work over
$\ZZ$  we would be required to use $2\alpha: \ll(V) \to \ll(V)$ since
in the free Lie algebra over $\ZZ$ $[v_i,v_i] = 2 v_i \otimes v_i$.

To simplify notation, identify the tensor algebra with the free
non-commutative algebra, so we write $vw$ and $v^2$ in place of $v
\otimes w$ and $v \otimes v$ respectively.

Let us compute $d$ on basis elements. By \eqref{geoderiv} and
\eqref{derivation}, we see that
    \begin{equation}\label{computationofd}
        \begin{split}
        d(v_{i_1} v_{i_2} \ldots v_{i_n})   &= d(v_{i_1}) v_{i_2} \ldots v_{i_n}  - v_{i_1} d(v_{i_2} v_{i_3} \ldots
        v_{i_n})\\
                                &= d(v_{i_1}) v_{i_2} \ldots v_{i_n} -
                                v_{i_1}
                                d(v_{i_2}) v_{i_3} \ldots v_{i_n} + v_{i_1}
                                v_{i_2}
                                d(v_{i_3} \ldots v_{i_n}) \\
                                &= d(v_{i_1}) v_{i_2} \ldots v_{i_n} -
                                v_{i_1}
                                d(v_{i_2}) v_{i_3} \ldots v_{i_n} + \ldots
                                + \\ &
                                (-1)^{(n-1)} v_{i_1} v_{i_2} \ldots v_{i_{n-1}} d(
                                v_{i_n})\\
                                &= v_{i_1}^2 v_{i_2} \ldots v_{i_n} -
                                v_{i_1}
                                v_{i_2}^2 v_{i_3} \ldots v_{i_n} + \ldots \\ & +
                                (-1)^{(n-1)} v_{i_1} v_{i_2} \ldots
                                v_{i_{n-1}}
                                v_{i_n}^2
        \end{split}
    \end{equation}

Notice that by construction, the map of \eqref{computationofd} is a
derivation of Lie algebras. Thus $d: \ll(V) \to \ll(V)$ is a Lie
algebra derivation of $\ll(V)$ of degree one  and $D: \ll(V) \to
\ll(V)$ is a Lie algebra derivation of $\ll(V)$ of degree zero.

The free Lie algebra $\ll(V)$ has bracket
    \[ [,]: L_p(V) \otimes L_q(V) \to L_{p+q}(V) \]
which satisfies the usual graded Lie algebra relations:
    \begin{enumerate}
            \item (graded antisymmetry)
            \begin{equation}
                [x,y]= (-1)^{|x||y|+1}[y,x]
            \end{equation}
            \item (graded Jacobi)
            \begin{equation}
                [x,[y,z]]=[[x,y],z]+(-1)^{|x||y|}[y,[x,z]]
            \end{equation}
            or equivalently
            \begin{equation}
            (-1)^{|x||z|}[x,[y,z]] + (-1)^{|z||y|)}[z,[x,y]]+(-1)^{|x||y|}[y,[z,x]]=0
            \end{equation}
    \end{enumerate}
On elements of the form $v_{i_1} \tensor \ldots \tensor v_{i_p} \in
V^{\tensor p}$, $v_{j_1} \tensor \ldots \tensor v_{j_q} \in
V^{\tensor q}$, it is a computational exercise to show that the
bracket
  \begin{equation}
        [,]: V^{\tensor p} \otimes V^{\tensor q}   \to V^{\tensor
        p+q},
  \end{equation}
  \begin{multline}
        (v_{i_1} \tensor \ldots \tensor v_{i_p}) \otimes (v_{j_1} \tensor \ldots \tensor v_{j_q})  \mapsto
         v_{i_1} \tensor \ldots \tensor v_{i_p} \otimes v_{j_1} \tensor \ldots \tensor
        v_{j_q} \\
        -(-1)^{pq} v_{j_1} \tensor \ldots \tensor v_{j_q} \tensor v_{i_1} \tensor \ldots \tensor v_{i_p}
  \end{multline}
may be expressed in terms of the group ring $\ZZ[\Sigma_{p+q}]$ as
right multiplication by
    \begin{equation}\label{defB}
        B_{p,q}=(1 - (-1)^{pq} \sigma_{(p+q, p+q-1, \ldots,
        2,1)}^{p}).
    \end{equation}

 Suppose $a \in L_p$ and $b \in L_q$. Then $a$ (respectively $b$) is a linear combination of elements
 of the form $R_{w_p}(v_{i_1} \tensor \ldots \tensor
 v_{i_p})$ (respectively $R_{w_q}(v_{j_1} \tensor \ldots \tensor v_{j_q})$. Since
 the bracket maps $L_p \tensor L_q$ into $L_{p+q}$, we have that
    \begin{equation}\label{BisinL}
        R_{B_{p,q}} (R_{w_p}(v_{i_1} \tensor \ldots \tensor v_{i_p}) \tensor R_{w_q}(v_{j_1} \tensor
        \ldots \tensor v_{j_q})) \in L_{p+q},
    \end{equation}
which means the expression of \eqref{BisinL} can be written as a
linear combination of elements of the form $R_{w_{p+q}}(v_{k_1}
\tensor \ldots \tensor v_{k_{p+q}})$.

Summarizing, we have seen that if $S$ is interpreted as a basis
concentrated in degree one, then for each $n,p,q \in \NN$ we have the
following maps:
\begin{enumerate}
    \item
    \begin{equation}\label{eq:paul}
        R_{w_n}: \Hom_{\dd_\fin}([n],S) \to \Hom_{\dd_\fin}([n],S)
    \end{equation}
    which yielded the derivation of degree $0$.
    \item
    \begin{equation}\label{eq:sasha}
        d: \Hom_{\dd_\fin}([n],S) \to \Hom_{\dd_\fin}([n+1],S)
    \end{equation}
    which gave a derivation of degree $1$.
    \item
    \begin{equation}\label{eq:john}
        R_{B_{p,q}}: \Hom_{\dd_\fin}([p]\amalg [q],S) \to \Hom_{\dd_\fin}([p]\amalg[q],S)
    \end{equation}
        which, after identifying $\Hom_{\dd_\fin}([p],S) \tensor
        \Hom_{\dd_\fin}([q],S)$ with $\Hom_{\dd_\fin}([p]\amalg
        [q],S)$ gave the bracket of the free Lie algebra on
        $\ZZ[S]=\Hom_{\dd_\fin}([1],S)$.
\end{enumerate}

\section{Construction of the d.g.l.}

This section will develop the identities, expressed as commutative
diagrams, in $\dd = \dd_\fin$ which allow us to construct a d.g.l.

\begin{defn}
    A {\bf complex} in $\dd$ is a collection of objects and
    morphisms in $\dd$, $(A_n,d_n)$, indexed by $\ZZ$, such that
    $d_n \in \Hom_\dd(A_n, A_{n-1})$ and $d_n \circ d_{n-1}=0 \in \Hom_D(A_n,
    A_{n-2})$. The integer, $n$, corresponding to $A_n$ may be
    called the degree of $A_n$.
\end{defn}

Given $S, T \in \fin$ then $\Hom_\fin(T,S) \in \fin$. If $f \in
\Hom_\dd(T,T')$ then $f^*: \Hom_\fin(T',S) \to \Hom_\fin(T,S)$ is a
morphism in $\dd$ since it is a linear combination of set maps
between $\Hom_\fin(T,S)$ and $\Hom_\fin(T',S)$.

For $Y\in \dd$, define $\End_\dd(Y)= \Hom_\dd(Y,Y)$. We can define
the action of a group $G$ on an object $Y \in\dd$ as a ring
homomorphism $\ZZ[G] \to \End_\dd(Y)$. In our case $\Sigma_n$ acts on
$[n] \in \dd$, thus for any $S\in\fin$, $\Sigma_n$ acts on
$\Hom_\fin([n],S)$, so we have a homomorphism $\ZZ[\Sigma_n] \to
\End_\dd(\Hom_\fin([n],S))$. With this homomorphism, $w_n \in
\ZZ[\Sigma_n]$ gives an element of $\End_\dd( \Hom_\fin([n], S))$.

\subsection{The ``geometric'' differential on the level of finite
sets.}

First we develop the idea of how the geometric cobar complex may be
considered as arising from a diagram in $\dd$.

For each $i \in \{ 1 ,\ldots , n+1 \}$, let $\delta_i \in
\Hom_\fin([n+1], [n])$ be given by
    \begin{equation} \label{eq:simpid}
        \delta_i(j) = \begin{cases} j-1 & \text{ if } i \neq j \\
                                    j & \text{ if } i=j \end{cases}
    \end{equation}

 Let $f_n \in \Hom_\dd([n+1], [n])$ be the map given by
    \begin{equation} \label{eq:deffn}
     f_n = \sum_{i=1}^n (-1)^{i-1} \delta_i.
    \end{equation}
 It is the alternating sum of the $\delta_i$, and as a linear combination
 of maps in $\fin$, it is a map in $\dd$.

Let $G: \fin^\opp \to \cc$ be any of the four functors, $S(\bullet)$,
$\HomS$, $X^\bullet$, or $F(\bullet)$ from $\fin^\opp$ to $\fin$,
vector spaces over $k$, pairs, or $\aa$. We consider each case of
applying $G$ to the complex in $\dd$:
    \begin{equation}\label{simpcobar}
    [1] \xleftarrow{f_1} [2] \xleftarrow{f_2} [3] \xleftarrow{f_3} \ldots
    \end{equation}
and for fixed $N \in \NN$, the complex
    \begin{equation}\label{simpcobarN}
    [1] \xleftarrow{f_1} [2] \xleftarrow{f_2} [3] \xleftarrow{f_3}
    \ldots\xleftarrow{f_{N-1}}
    [N].
    \end{equation}

Now $G$ can be applied to \eqref{eq:deffn}. Using simplicial
identities which arise from \eqref{eq:simpid}, one can show directly
that $f_{n+1}\circ f_n=0$. Alternatively, the following introduces
the technique of proof that will follow.

\begin{prop} $f_n \circ f_{n+1}= 0$ in the category $\dd$.
\end{prop}
\pf
    For $S \in \fin$, apply the functor $\Hom_{\dd}(\bullet,
    S)$. Then
    \begin{equation}
        \Hom_{\dd}( [n], S) \xrightarrow{ \Hom_{\dd}( f_n,
        S)}
        \Hom_{\dd}( [n+1], S) \xrightarrow{ \Hom_{\dd}( f_{n+1},
        S)}
        \Hom_{\dd}( [n+2], S).
    \end{equation}
    After identifying $\Hom_{\dd}( [n], S)$ with
    $\ZZ[S]^{\tensor n}$ we have
    \begin{equation}
        \ZZ[S]^{\tensor n} \xrightarrow{d_n}
        \ZZ[S]^{\tensor n+1} \xrightarrow{d_{n+1}}
        \ZZ[S]^{\tensor n+2}
    \end{equation}
    where $d_n$ is the algebra derivation of \eqref{computationofd}
    (where $\oplus_n\ZZ[S]^{\tensor n}$ is the algebra). Then the
    computation of $d \circ d$ on generators $s_i \in S$,
    \begin{equation}
        d \circ d (s_i) = d ( s_i^2) = (ds_i) s_i + (-1)^{1\cdot
        1}s_i \cdot (d s_i) = s_i^3 - s_i^3=0
    \end{equation}
    shows that $d^2=0$.

    Since this holds for all $S \in \fin$, $f_n \circ f_{n+1}$ gives a natural
    transformation of functors
    \begin{equation} \label{eq:bomb}
        \Hom_{\dd}( [n], \bullet) \to \Hom_{\dd}( [n+2], \bullet)
    \end{equation}
    which is now seen to be identically $0$.
    Now applying \eqref{eq:bomb} to the object $[n] \in \fin$, and
    taking the image of the identity morphism in $\Hom_{\dd}( [n],
    [n])$ we see that $f_n \circ f_{n+1} =0$.
    \endpf

In the above proof we are writing out in detail the proof of a
special case of the Corollary to Yoneda's Lemma (see \cite{mac}
Chapter 3), that every natural transformation between two
representable functors $\Hom_\cc(X,\bullet)$ and
$\Hom_\cc(Y,\bullet)$ is of the form $\Hom_\cc(h, \bullet)$ for some
$h: Y \to X$. Thus since the natural transformation in the proof is
given by $\Hom_\dd(f_n \circ f_{n+1}, \bullet)$, and the functor is
zero, we conclude that $f_n \circ f_{n+1} = 0$.

In the case where we take $G=X^\bullet$ we have a version of the
geometric cobar construction. Since $f_{n+1} \circ f_{n} = 0$, we
have $G(f_n) \circ G(f_{n+1})=0$, and we have a complex in each
category. The case of primary interest is of course when $G=F$ and we
have
    \begin{equation}\label{eq:star}
        F([1]) \xrightarrow{Ff_1} F([2]) \xrightarrow{Ff_2} \ldots
    \end{equation}
    and
    \begin{equation}\label{eq:starN}
        F([1]) \xrightarrow{Ff_1} F([2]) \xrightarrow{Ff_2} \ldots
        \xrightarrow{Ff_{N-1}} F([N])
    \end{equation}
which are objects of $\Ch(\aa)$, chain complexes in $\aa$. We will
abbreviate $F(n)$ for $F([n])$.

 We may assemble \eqref{eq:star} and define $R_F = \oplus_{n=1}^\infty
F(n)$ and from \eqref{eq:starN} define $R_F^N = \oplus_{n=1}^N F(n)$.
(We abbreviate $R_F$ and $R_F^N$ by $R$ and $R^N$ when $F$ is
understood.) The tensor product in $\aa$ gives each of these a
product structure, where in the case of \eqref{eq:starN}, if $p+q>N$
then $F(p) \times F(q) \xrightarrow{\tensor} F(p+q) \equiv 0$ is
defined to be zero.

The tensor product structure is associative by axiom
\eqref{associativity},  thus making $R$ and $R^N$ associative
algebras.

Now to see that the tensor product makes \eqref{eq:star} and
\eqref{eq:starN} differential graded algebras,  we must show that
    \begin{equation*}
        R \times R \xrightarrow{\tensor} R
        \quad (\text{or } R^N \times R^N \xrightarrow{\tensor} R^N)
    \end{equation*}
preserves the differential. That is for each $p,q \in \NN$,
    \begin{equation}\label{eq:wisemagic}
    \xymatrix{
    F(p)\tensor F(q) \ar[r] \ar[d] & F(p+q)\ar[d] \\
    F(p+1)\tensor F(q) \oplus F(p) \tensor F(q+1) \ar[r] & F(p+q+1)
    }
    \end{equation}
commutes, where the top map is given by $N(p,q)$, the left map by
$(Ff_p \tensor 1 , (-1)^p 1 \tensor Ff_q)$, the right map by
$Ff_{p+q}$, and the bottom map by $N(p+1,q) + N(p,q+1)$.

\begin{prop} $R$ (or $R^N$) is a d.g.a.
\end{prop}
    \pf We must show that the two maps of \eqref{eq:wisemagic} are
    equal. $F(p) \tensor F(q) \in \aa$ is the functor $\tensor \circ F \times
    F$ applied to $[p] \times [q] \in \fin \times \fin$. To obtain
    our top right map we apply $N(p,q)$ and the map $F(f_{p+q})$:
    \begin{equation} \label{eq:slash}
        F(p) \tensor F(q) \xrightarrow{N} F( p \amalg q)
        \xrightarrow{F(f_{p+q})} F( p + q+ 1)
    \end{equation}

    The following commutes by definition of natural transformation.
      \begin{equation}  \xymatrix{
            F(p) \tensor F(q) \ar[r]^N \ar[d]^{F(f_p)\tensor 1} & F(
            p \amalg q) \ar[d]^{F(f_p \times \id_q)} \\
            F(p+1) \tensor F(q) \ar[r]^N & F(p+q+1)
        }
      \end{equation}

    Therefore we may replace the maps in \eqref{eq:wisemagic} going around
    the bottom left by the sum of the two maps
    \begin{gather} \label{eq:axel}
        F(p) \tensor F(q) \xrightarrow{N} F( p \amalg q)
        \xrightarrow{F(f_{p} \amalg \,  \id_q)} F( p + q+ 1) \\
        \label{eq:izzy}
        F(p) \tensor F(q) \xrightarrow{N} F( p \amalg q)
        \xrightarrow{F( (-1)^q \id_p \amalg f_{q})} F( p + q+ 1)
    \end{gather}

    We need to prove that $(F( f_p \amalg \Id_{q}) + F( (-1)^q \id_p \amalg
    f_{q})) \circ N = F(f_{p+q}) \circ N$. This will follow if we
    show that
    \begin{equation}\label{eq:luminaire}
    f_p \amalg \, \Id_q + (-1)^p \Id_p \amalg f_q = f_{p+q} \in
    \Hom_{\dd}([p+q+1],[p+q]).
    \end{equation}

    To prove \eqref{eq:luminaire}, for fixed $p, q \in \NN$, consider the
    two representable bi-functors, $\tensor \circ (\Hom_\dd([p], S)
    \times \Hom_\dd([q], T))$ and $\Hom_\dd([p+q+1], S \amalg T)$, for
    $S,T \in \fin$. These are functors from $\fin \times \fin $ to
    $\ZZ$--modules.
            Furthermore,  $\tensor \circ (\Hom_\dd([p], S)
            \times \Hom_\dd([q], T))$ can be identified with the functor
            $\Hom_\dd([p] \amalg [q], S \amalg T)$ (since they are
            identified as Abelian groups).
    Consider the natural transformation of functors given by
    $\Hom_\dd(f_p \amalg \, \Id_q + (-1)^p \Id_p \amalg f_q - f_{p+q},
    \bullet)$. We will show that this transformation is the zero map
    for all $S$ and $T$.  Taking $S=[p]$, $T=[q]$ with $\Id_p \amalg
    \Id_q \in \Hom_\dd([p]\amalg [q],[p]\amalg [q])$ proves \eqref{eq:luminaire}.

    After identifying the functor $\Hom_{\dd}([n],S)$ with
    $\ZZ[S]^{\tensor n}$, to show that the natural transformation
    above is zero amounts to checking that, for all $S,T \in \fin$, the
    diagram
    \begin{equation}
        \xymatrix{
        \ZZ[S]^{\tensor p} \tensor \ZZ[T]^{\tensor q} \ar[r] \ar[d] & \ZZ[S
        \amalg T]^{\tensor p+q} \ar[d] \\
        \ZZ[S]^{\tensor p+1} \tensor \ZZ[T]^{\tensor q} \oplus
        \ZZ[S]^{\tensor p} \tensor \ZZ[T]^{\tensor q+1} \ar[r] &
        \ZZ[S \amalg T]^{ \tensor p+q+1}
        }
    \end{equation}
    commutes. This is a straightforward calculation using the
    definition of the derivation of \eqref{computationofd}.
    \endpf

As differential graded algebras $R$ and $R^N$ acquire the structure
of differential graded Lie algebras in their own right, just as any
associative algebra may be considered as a Lie algebra. Now $R$ is
the universal enveloping algebra of a primitive Lie algebra, which we
will construct below.

We wish to produce a Lie subalgebra of $R$, and will make use of our
work in Section \ref{sec:algident} to do so.

We saw that for any $n\in \NN$, $S \in \fin$ we have a map
    \begin{equation}
        R_{w_n}: \Hom_{\dd}([n],S) \to \Hom_{\dd}([n],S)
    \end{equation}
    which is given by right multiplication by $w_n$. As done before,
    this can be seen as a natural transformation of representable
    functors, and so $R_{w_n}$ is represented as
    $\Hom_{\dd}(f,S)$ for some morphism $f$ in $\dd$.
    Denote this morphism by $\widetilde{w_n}$. Again one may wish to
    think of this as the image of the identity map in the case of
    $S=[n]$.

    So $\widetilde{w_n}: [n] \to [n]$, and as a morphism in $\dd$,
    applying the functor $\Hom_\dd([n], \bullet)$, we
    get the map $R_{w_n}$.
    Furthermore, from Proposition \ref{wn2wn}, we can conclude that
    $\widetilde{w_n}^2 = n \widetilde{w_n}$.

Now take  $G=F$ and $\cc= \aa$. At this point we assume that $\aa$
has the following property:
\begin{defn}
An additive category $\aa$ is $\QQ$-Karoubian if, for any morphism
$f: M \to M$ in $\aa$ such that there is an integer $n$ such that $f
\circ f = n \cdot f$  the image of $f$ is an object in $\aa$.
\end{defn}

Notice that if $\aa$ is Abelian this condition is satisfied. If $\aa$
is Karoubian then $f^2=f$ implies the image of $f$ is in $\aa$, so
the Karoubian requirement is not quite enough.

For each $n$ we have the right multiplication by $F(\widetilde{w_n})$
on $F([n])$. By Proposition \ref{wn2wn},  the $\QQ$-Karoubian
hypothesis implies that the image of $F(\widetilde{w_n})$ is an
object of $\aa$. For the functor $\Hom_\dd([n], \bullet)$ there is
also no problem in taking images, though we will need to invert $2$
to get all our desired properties. However, when working with the
functor $X^\bullet$, we cannot take the image of
$X(\widetilde{w_n})$.

The rest of this section is dedicated to proving the following
theorem:

\begin{thm}\label{thm:maindgl} Let $\pp_F = \oplus_{n=1}^\infty F(n)F(\widetilde{w_n})$ and let
$\pp_F^N = \oplus_{n=1}^N F(n)F(\widetilde{w_n})$. $\pp_F$ and
$\pp_F^N$ are differential graded Lie algebras in $\aa$.
\end{thm}

\subsection{Closure under the differential}  We would first like to show
$\pp_F$ and $\pp_F^N$ are closed under the differential $F(f_n)$:

To show that $\pp_F$ and $\pp_F^N$ are subcomplexes we must show that
we have a map
    \begin{equation}
        F(f_n): F(n)F(\widetilde{w_n}) \to F(n+1)F(\widetilde{w_{n+1}}).
    \end{equation}
That is, since we already know $F(f_n) \circ F(f_{n+1})=0$, we are
required to show that
    \begin{equation}\label{eq:subcomplexcondition}
        F(f_n) ( F(n)F(\widetilde{w_n})) \subset F(n+1)F(\widetilde{w_{n+1}}).
    \end{equation}

If we show that there exists $\phi_n$ such that
    \begin{equation}\label{eq:digweed}
        \xymatrix{
        [n+1] \ar[r]^{f_n} \ar[d]_{\widetilde{w_{n+1}}} & [n]
        \ar[d]^{\widetilde{w_n}} \\
        [n+1] \ar[r]^{\phi_n} & [n]
        }
    \end{equation}
    commutes in $\dd$,
then applying $F$, we have a commutative diagram
    \begin{equation}
        \xymatrix{
            F(n) \ar[r]^{F(\phi_n)} \ar[d]_{F(\widetilde{w_{n}})} &
            F(n+1) \ar[d]^{F(\widetilde{w_{n+1}})} \\
            F(n) \ar[r]^{F(f_n)} & F(n+1)
            }
    \end{equation}
which, in particular, shows that \eqref{eq:subcomplexcondition}
holds. It is important to note that $\phi_n$ is {\it not} the same as
$f_n$.

The following proposition is the first step in showing
\eqref{eq:digweed} commutes.

\begin{prop} \label{prop:mainpropfin2}
     There exists $g_n \in \dd$ such that the following diagram commutes for all $S \in \fin$
    \begin{equation}\label{simptougheqn}
            \xymatrix{
            \Hom_{\dd}([n],S)  \ar[r]^{g_n} \ar[d]^{w_{n}} & \Hom_{\dd}([n+1], S)
            \ar[d]^{w_{n+1}} \\
            \Hom_{\dd}([n],S) \ar[r]^{f_{n}^*} &\Hom_{\dd}([n+1],S)
            }
    \end{equation}
\end{prop}
\begin{remark}
    In the course of the proof we will see that $2$ needs to be
    inverted to make this Proposition hold.
\end{remark}
\pf
    Suppose $S = \{ s_1, \ldots, s_m \}$.
    \eqref{simptougheqn} may be identified with
    \begin{equation}
            \xymatrix{
            \ZZ[S]^{\tensor n}  \ar[r]^{g_n} \ar[d]^{w_{n}} &
            \ZZ[S]^{\tensor n+1}
            \ar[d]^{w_{n+1}} \\
            \ZZ[S]^{\tensor n} \ar[r]^{d} & \ZZ[S]^{\tensor
            n+1}
            }
    \end{equation} where we must solve for $g_n$.

    Let us look at $ {2}d$ on generators of $\ll( \ZZ[S])$:
    $ {2}d(s_i) = [s_i,s_i]$, and thus \[ {2}d
    (\ll(\ZZ[S])) \subset \ll(\ZZ[S]). \] Then
    \begin{gather}
    \frac{1}{2}d(R_{w_n}(s_{j_1}\tensor s_{j_2} \tensor \ldots \tensor
    s_{j_n})) \\
    = \frac{1}{2}d ( [s_{j_1} [ s_{{j_2}} [ \ldots [ s_{j_{n-1}},
    s_{j_n}] \ldots ]) \in L_{n+1}(\ZZ[S]) \\
    = R_{w_{n+1}}( \text{ an expression in $s_j$'s of degree }n+1) \\
    = R_{w_{n+1}}(g_n( s_{j_1}\tensor s_{{j_2}} \tensor \ldots \tensor
    s_{j_n}))
    \end{gather}
    which defines $g_n$.
\endpf

Since Proposition \ref{prop:mainpropfin2} holds for all $S$, taking
$S=[n]$ and looking at the image of the identity map under $g_n$, we
obtain $\phi_n \in \Hom([n+1],[n])$. Thus we have \eqref{eq:digweed},
completing the proof that $\pp_F$ and $\pp_F^N$ are complexes.

~

\subsection{Closure under the bracket}

Now to construct a Lie algebra we need to show that taking the images
defined by the $F(\widetilde{w_n})$ will be closed under the bracket.

Recall from \eqref{eq:john} that right multiplication by $B_{p,q}$,
$R_{B_{p,q}}: \Hom_\dd( [p]\amalg[q], S) \to \Hom_\dd( [p]\amalg[q],
S)$ gave the bracket on $\ll(\ZZ[S])$. Considering this as a natural
transformation of functors, we again can conclude that there exists
$\widetilde{B_{p,q}} \in \Hom_\dd([p] \amalg [q],[p] \amalg [q])$,
such that $B_{p,q} = \Hom_\dd( \widetilde{B_{p,q}}, S)$.

Since $R$ has an algebra structure, we can define a bracket by taking
the multiplication map
    \begin{equation}
        F(p) \times F(q) \xrightarrow{\tensor} F(p)\tensor F(q)
        \xrightarrow{N(p,q)} F(p+q)
    \end{equation}
and adding the map
    \begin{equation}
        F(p) \times F(q) \xrightarrow{\tensor} F(p)\tensor F(q)
        \xrightarrow{-(-1)^{pq}\sigma} F(q) \tensor F(p)
        \xrightarrow{N(q,p)} F(p+q).
    \end{equation}
Note that as a map  $[,]: F(p) \tensor F(q) \to F(p+q)$, $[,] =
N(p,q) - (-1)^{pq}\sigma N(q,p)$. Our goal is to show that this map
sends $F(p)F(\widetilde{w_p}) \tensor F(q)F(\widetilde{w_q})$ into
$F(p+q)F(\widetilde{w_{p+q}})$.

By the definition of $\widetilde{B_{p,q}}$ \eqref{defB} and the
commutativity axiom condition on $F$, we have that the following
diagram commutes:
    \begin{equation}\label{eq:silver}
        \xymatrix{
            F(p) \tensor F(q) \ar[r]^{N(p,q)} \ar[dr]_{[,]} & F(p+q)
            \ar[d]^{F(\widetilde{B_{p,q}})}   \\
            & F(p+q)
        .}
    \end{equation}

Since $N$ is a natural transformation, the following diagram
commutes:
    \begin{equation}\label{eq:copper}
    \xymatrix{
        F(p) \tensor F(q) \ar[r]^N \ar[d]_{F( \widetilde{w_p} ) \tensor F( \widetilde{w_q})} &
        F( [p] \amalg [q]) \ar[d]^{F( \widetilde{w_p} \amalg  \widetilde{w_q})}   \\
        F(p) \tensor F(q) \ar[r]^N & F( [p] \amalg [q])
        ,}
    \end{equation}
thus the image of $F(p) F(  \widetilde{w_p} ) \tensor F(q) F(
\widetilde{w_q}) $ under $N$ is contained in the image of $F(
\widetilde{w_p} \amalg  \widetilde{w_q})$.

Thus if we show that for each $p,q \in \NN$ there exists a map
$\widetilde{\psi_{p,q}}$ so that the following diagram in $\dd$
commutes:
    \begin{equation}\label{eq:gatecrasher}
        \xymatrix{
        [p] \amalg [q] & &  [p] \amalg [q] \ar[ll]_{\widetilde{\psi_{p,q}}} \\
        [p] \amalg [q] \ar[u]^{\widetilde{w_p} \amalg \widetilde{w_q}}
         & &  [p] \amalg [q] \ar[ll]_{\widetilde{B_{p,q}}}
         \ar[u]_{\widetilde{w_{p+q}}}
        ,}
    \end{equation}
then applying $F$ yields
    \begin{equation}\label{eq:gold}
        \xymatrix{
        F([p] \amalg [q]) \ar[d]_{F(\widetilde{w_p} \amalg \widetilde{w_q})}
        \ar[rr]^{F(\widetilde{\psi_{p,q}})} & &  F([p] \amalg [q]) \ar[d]^{F(\widetilde{w_{p+q}})}  \\
        F([p] \amalg [q]) \ar[rr]^{F(\widetilde{B_{p,q}})}
         & &  F([p] \amalg [q])
        ,}
    \end{equation}
which is enough to show that the bracket maps $F(p)F(\widetilde{w_p})
\tensor F(q)F(\widetilde{w_q})$ into $F(p+q)F(\widetilde{w_{p+q}})$.

Once we have Proposition \ref{prop:finbracket} below, we prove
\eqref{eq:gatecrasher} by arguing as before, taking $S=[p]\amalg [q]$
and looking at the image of the identity map in $\Hom_\dd([p]\amalg
[q], [p]\amalg[q])$ under $\widetilde{\psi_{p,q}}$.

\begin{prop} \label{prop:finbracket}
    There exists $\widetilde{\psi_{p,q}} \in \dd$ such that for any $S\in \fin$
    \begin{equation} \label{eq:finbracket}
        \xymatrix{
        \Hom_\dd([p] \amalg [q], S) \ar[d]_{\Hom_\dd(\widetilde{w_p} \amalg \widetilde{w_q},S)}
        \ar[rr]^{\Hom_\dd(\widetilde{\psi_{p,q}},S)} & &  \Hom_\dd([p] \amalg [q],S) \ar[d]^{\Hom_\dd(\widetilde{w_{p+q}},S)}  \\
        \Hom_\dd([p] \amalg [q],S) \ar[rr]^{\Hom_\dd(\widetilde{B_{p,q}},S)}
         & &  \Hom_\dd([p] \amalg [q],S)
        }
    \end{equation}
    commutes.
\end{prop}
\pf Suppose $S = \{ s_1, \ldots, s_m \}$.
    \eqref{eq:finbracket} may be identified with the following diagram where we are trying to solve for
    $\psi_{p,q}$ such that
    \begin{equation}
        \xymatrix{
            \ZZ[S]^{\tensor p} \tensor \ZZ[S]^{\tensor q}
            \ar[rr]^{\psi_{p,q}} \ar[d]_{R_{w_p} \tensor R_{w_q}} & &
            \ZZ[S]^{\tensor p+q} \ar[d]^{R_{w_{p+q}}} \\
            \ZZ[S]^{\tensor p} \tensor \ZZ[S]^{\tensor q}
            \ar[rr]^{[,]} & & \ZZ[S]^{\tensor p+q}
        }
    \end{equation}
    commutes.
    Now the image of the left map is $L_p(\ZZ[S]) \tensor
    L_q(\ZZ[S])$, and the bracket maps this to $L_{p+q}(\ZZ[S])$.
    Thus given generators $s_{i_1}\tensor s_{i_2} \tensor \ldots \tensor
    s_{i_p}$ and $s_{j_1}\tensor s_{j_2} \tensor \ldots \tensor
    s_{j_q}$ of $\ZZ[S]^{\tensor p}$ and $\ZZ[S]^{\tensor q}$, the
    composite of the left and bottom map yields
    \begin{gather}
        [ R_{w_p}( s_{i_1}\tensor s_{i_2} \tensor \ldots \tensor
        s_{i_p}), R_{w_q}(s_{j_1}\tensor s_{j_2} \tensor \ldots \tensor
        s_{j_q})] \in L_{p+q}(\ZZ[S]) \\
        = R_{w_{p+q}}(\text{ an expression in the $s_i, s_j$}) \\
        = R_{w_{p+q}}( \psi_{p,q})
    \end{gather}
which defines $\psi_{p,q}$.
\endpf

Notice that since the differential is compatible with the algebra
structure (i.e. it is a graded algebra derivation), the differential
is automatically a graded Lie derivation, and thus $\pp_F$ and
$\pp_F^N$ are differential graded Lie algebras, thus completing the
proof of Theorem \ref{thm:maindgl}.

\begin{remark}
    $\pp^N_F$ is a Lie-quotient of $\pp_F$, as the kernel
    is the Lie ideal generated by all ``words'' of bracket length greater than
    $N$ (which is clearly stable under the bracket).
\end{remark}

\section{Applying the Construction to the Abelian Category of Motives}

Now we specialize to the category of motives. Here we summarize some
of the main facts concerning Nori's work:

\subsection{Summary of EHM and $\CLD$}

Fix an embedding of $k \hookrightarrow \CC$. All varieties and
morphisms discussed below will be defined over $k$. Let $\Zmod$
denote the category of all finitely generated Abelian groups, and let
$\Ab$ denote the category of all Abelian groups. All these
constructions are made with integer coefficients.

\begin{fact}
$\ehm$ is an Abelian category equipped with a faithful exact functor,
$\forg: \ehm \to \Zmod$.
\end{fact}

\begin{fact}
 For all varieties $X$ and closed subvarieties $Y \subset X$ (defined
 over $k$), and all non-negative integers $q$, there is an object
 $H_q(X,Y)$ of $\ehm$.
\end{fact}

\begin{fact}\label{fact:three}
    With pairs $(X,Y)$ and $(X',Y')$ as above, every commutative
    diagram
    \begin{equation*}
        \xymatrix{
            Y \ar[d] \ar[r]^g & Y' \ar[d] \\
            X \ar[r]^f & X' }
    \end{equation*}
    induces a morphism $H_q(f,g): H_q(X,Y) \to H_q(X',Y')$ in $\ehm$.
\end{fact}

\begin{fact} \label{fact:four}
    For $Z \subset Y \subset X$ closed subvarieties of $X$, there are
    boundary operators
    \begin{equation*}
        H_q(X,Y) \to H_{q-1}(Y,Z).
    \end{equation*}
\end{fact}

\begin{fact}
    Let $H_q(X(\CC),Y(\CC))$ denote the $q$-th singular homology of
    the pair \\ $(X(\CC), Y(\CC))$ with $\ZZ$ coefficients. Then
        \begin{equation*}
            \forg H_q(X,Y) \cong H_q(X(\CC),Y(\CC))
        \end{equation*}
\end{fact}

\begin{fact}
    Applying $\forg$ to Fact \ref{fact:three} yields the arrows on
    singular homology induced by $(f,g)$ and applying $\forg$ to Fact
    \ref{fact:four} yields the usual boundary operator.
\end{fact}

Because $\forg$ is faithful exact, from the standard properties of
singular homology we deduce Facts \ref{fact:seven} --
\ref{fact:nine}:

\begin{fact}\label{fact:seven}
        For $X \supset Y \supset Z \supset W$ (defined over $k$)
        and any $q \in \ZZ$, the composite of the boundary operators
        in $\ehm$
        \begin{equation*}
            H_q(X,Y) \to H_{q-1}(Y,Z) \to H_{q-2}(Z,W)
        \end{equation*}
        is zero.
\end{fact}

\begin{fact}
    There is the usual exact sequence of a triple $(X,Y,Z)$ in
    $\ehm$ (for the usual exact sequence of a triple in the topological setting see chapter 4 section 8 of \cite{sp}).
\end{fact}

\begin{fact}\label{fact:nine}
    With notation as in Fact \ref{fact:three} whenever the composite
    is defined we have
        \begin{equation*}
            H_q(f,g) \circ H_q(f',g') = H_q(f \circ f', g \circ g').
        \end{equation*}
\end{fact}

\begin{fact}\label{fact:eleven}
    There is a functor $\tensor : \ehm \times \ehm \to \ehm$ such
    that
for all objects $A,B$ in $\ehm$ we have
        \begin{equation*}
            \forg ( A \tensor B) = \forg A \tensor \forg B.
        \end{equation*}
This implies that
    the $\tensor$ structure on $\ehm$ satisfies the usual
    associativity and commutativity constraints, and
    $H_0(\spec(k),\emptyset) = \mathbf{1}$ serves as an identity.
\end{fact}

\begin{defn}
    An $r$-admissible pair consists of an affine variety $X$ and a
    closed subvariety $Y$ of $X$ such that $\dim Y < r$, $X - Y$ is
    smooth of pure dimension $r$, and
        \begin{equation*}
            H_s(X(\CC),Y(\CC)) = \begin{cases} 0 & \text{if } s\neq r
            \\                                  \text{a free Abelian
            group} & \text{if } s = r. \end{cases}
        \end{equation*}
\end{defn}

\begin{fact} \label{fact:thirteen}
    For $r_i$-admissible pairs $(X_i,Y_i)$ for $i= 1,2$ we have
        \begin{equation*}
            H_{r_1}(X_1,Y_1) \tensor H_{r_2}(X_2,Y_2) =
            H_{r_1+r_2}((X_1,Y_1) \times (X_2,Y_2)).
        \end{equation*}
\end{fact}

Now $\indehm$ is the category obtained by taking direct limits of
objects in $\ehm$. (For a summary of ``Ind'' categories see
\cite{deligne}.) The properties of $\forg$ imply that the objects of
$\ehm$ satisfy the ascending chain condition. It follows that
$\indehm$ is an Abelian category. We denote by $\forg$, once again,
the induced functor $\forg: \indehm \to \Ab$.

The form of the lemma below is taken from \cite{nori}. It was first
proved by Beilinson \cite{be}.
\begin{thm}[The Basic Lemma] Let $X$ be an affine variety of
dimension $n$, and let $W \subset X$ be a closed subvariety of
dimension less than $n$. Then there exists a closed subvariety $Z$,
such that $X \supset Z \supset W$ and $(X,Z)$ is an $n$-admissible
pair.
\end{thm}

Let $X$ be an affine variety of dimension $n$. An admissible
filtration $\ff$  of $X$ is a filtration of $X$ by Zariski closed
sets $X = X_n \supset X_{n-1} \supset \ldots \supset X_0 \supset
X_{-1}= \emptyset$ such that $(X_i, X_{i-1})$ is an $i$-admissible
pair. An admissible filtration $\ff$ for $X$ defines a complex
    \begin{equation*}
        H_n(X,X_{n-1}) \to H_{n-1}(X_{n-1},X_{n-2}) \to \ldots \to
        H_0( X_1, \emptyset)
    \end{equation*}
    which we call $\CLD(X,\ff)$. Observe that $H_q(\forg \CLD(X,\ff)) =
    H_q(X(\CC))$ for the same reason that cellular homology agrees with
    singular homology for cell complexes. The Basic Lemma assures us
    that any filtration $\{ X_q \}_q$ with $\dim X_q \leq q$ is
    contained in an admissible filtration. Consequently admissible
    filtrations form a direct system. The direct limit of
    $\CLD(X,\ff)$ taken over all admissible filtrations is a complex
    in $\indehm$ called $\CLD(X)$.

    In fact $X \mapsto \CLD(X)$ is a functor from affine varieties
    over $k$ to chain complexes in $\indehm$. One may observe that a closed embedding $Y
    \hookrightarrow X$ gives a monomorphism $\CLD(Y) \to
    \CLD(X)$ and we may define $\CLD(X,Y) = \CLD(X) / \CLD(Y)$. Again
    the Basic Lemma is used to show that this is a functor since
    given $f: X \to Y$ for any filtration $\{ X_q \}$ applying $f$ and taking Zariski closures
    gives a filtration $\{\overline{f(X_q)} \}$ for which we can find an admissible
     filtration on $Y$ which dominates it.  Thus in the limit we have
     a map $\CLD(X) \to \CLD(Y)$.

    Given $\ff$ an admissible filtration on $X$ and $\gg$ an
    admissible filtration on $Y$, then $ \ff \times \gg$ is an
    admissible filtration on $X \times Y$, where
        \begin{equation*}
            (\ff \times \gg)_m = \bigcup_{p+q=m}\ff_pX \times \gg_qY.
        \end{equation*}

    It follows from Fact \ref{fact:thirteen} that
        \begin{equation*}
            \CLD(X,\ff) \tensor \CLD(Y,\gg) = \CLD(X\times Y, \ff
            \times \gg).
        \end{equation*}
    Taking direct limits over all $\ff$ and $\gg$ we get
        \begin{equation*}
            \CLD(X) \tensor \CLD(Y) \to \CLD( X \times Y)
        \end{equation*}
        in $\indehm$. The tensor structure gives a functor $\Ch
        (\indehm) \times \Ch(\indehm) \xrightarrow{\tensor}
        \Ch(\indehm)$ which satisfies the usual commutativity and
        associativity constraints of the tensor product -- in
        particular there is a  $\Sigma_n$-equivariant map
        $\CLD(X)^{\tensor n} \to \CLD(X^n)$ (refer to Proposition \ref{cortoaxioms}).
        That is, we have a natural transformation from the functor
        $\tensor \circ (\CLD \times \CLD): \Aff_k \times \Aff_k \to
        \Ch(\indehm)$ to the functor $\CLD \circ \times : \Aff_k
        \times \Aff_k \to \Ch(\indehm)$.

The lemma of homological algebra below yields subcomplexes
    \begin{equation*}
        \sing(X(\CC)) \supset P_*(X) \supset Q_*(X)
    \end{equation*}
such that both of the arrows
    \begin{equation*}
        \forg \CLD(X) = P_*(X)/Q_*(X) \twoheadleftarrow P_*(X)
        \hookrightarrow \sing(X(\CC))
    \end{equation*}
are quasi-isomorphisms.

Furthermore, both
    \begin{equation*}
        \xymatrix{
            P_*(X) \tensor P_*(Y) \ar@{^(->}[r] \ar[d] &
            \sing(X(\CC))\tensor \sing(Y(\CC)) \ar[d]^{EZ} \\
            P_*(X \times Y) \ar@{^(->}[r] & \sing(X(\CC) \times Y(\CC)) }
    \end{equation*}
and
    \begin{equation*}
        \xymatrix{
            \forg \CLD(X) \tensor \forg \CLD(Y) \ar[d] &
            P_*(X(\CC))\tensor P_*(Y(\CC)) \ar[d] \ar[l] \\
            \forg \CLD(X \times Y) & P_*(X(\CC) \times Y(\CC)) \ar[l]
            }
    \end{equation*}
commute.  Since $P_*, Q_*$ and $P_*/Q_*= \forg \CLD$ are all free
these maps are all quasi-isomorphisms.

Now $(\ehm, \forg)$ is universal for Facts 1 -- 6 and the tensor
structure is determined by Facts \ref{fact:eleven} and
\ref{fact:thirteen}. It follows that the d.g.l. constructed in
section 3 is a valid construction for $(\ehm, \forg)$. That is, by
taking the functor $F= \CLD(X^\bullet) \tensor \QQ$, $\pp_F$ is a
d.g.l. in $\mm$. Since we only use the properties above, the
construction is valid for mixed Hodge structures, Galois modules,
etc.

\begin{lemma}[Lemma of Homological Algebra]
    Suppose the following conditions are satisfied in an abelian
    category
    \begin{enumerate}
        \item $(A, \partial)$ is a chain complex
        \item $\{ F_pA \}_{p \in \ZZ}$ is an increasing filtration of
        subcomplexes of $A$
        \item for all $n \in \ZZ$, $\cup_p F_pA_n=A_n$
        \item for all $n \in \ZZ$ there exists $p \in \ZZ$ such that
        $F_pA_n = 0$
        \item for all $p \neq q$ in $\ZZ$, $H_q(gr^F_pA)=0$
    \end{enumerate}
    Then define $P_n$ to be the kernel of the composite map
    \begin{equation*}
        F_nA_n \to F_nA_{n-1} \to gr^F_nA_{n-1}
    \end{equation*}
    and set $Q_n = \partial (F_nA_{n+1}) + F_{n-1}A_n$.
    Then
        \begin{enumerate}
            \item $Q \subset P \subset A$ are subcomplexes.
            \item The quotient complex $P /Q$ is the familiar complex
            \begin{equation*}
                \ldots \to H_{n+1}(gr^F_{n+1}A) \to
                H_{n}(gr^F_{n}A) \to H_{n-1}(gr^F_{n-1}A)
            \end{equation*}
            \item $P \to P/Q$ and $P \to A$ are both
            quasi-isomorphisms.
        \end{enumerate}
\end{lemma}
This is a standard lemma, but as we were unable to find a reference
for it, we give the following proof due to Nori.

\pf The first two assertions are straightforward and rely only on
assumptions (i) and (ii). For the third we need the acyclicity of $Q$
and $A/P$. First recall that for a complex $(D,\partial)$ and an
integer $r$ we have the quotient complex $\tau_{<r}(D)$ given by
    \[ \ldots 0 \to D_r/Z_r(D) \to D_{r-1} \to D_{r-2} \to \ldots \]
    and the subcomplex $\tau_{>r}(D)$ given by
    \[ \ldots \to D_{r+2} \to D_{r+1} \to B_r(D) \to 0 \to 0 \ldots
    \]
    where $B_r(D) = \text{im}\partial (D_{r+1})$ and $Z_r(D)=
    \text{ker}(\partial: D_r\to D_{r-1})$.

    Now the increasing filtration $\{ F_pA \}_{p \in \ZZ}$ induces
    increasing filtrations $\{ F_p(A/P) \}$ and $\{ F_pQ \}$. We see that
    \[ gr^F_r(A/P)= \tau_{>r}(gr^F_rA) \] and
    \[ gr^F_r(Q)= \tau_{<r}(gr^F_rA). \]
The last hypothesis shows that the complexes $gr^F_r(A/P)$ and
$gr^F_r(Q)$ are acyclic. Then the five lemma shows that both
$F_aQ/F_bQ$ and $F_a(A/P)/F_b(A/P)$ are acyclic whenever $a>b$. Since
the third and fourth hypotheses hold for $Q$ and $A/P$ as well, we
deduce that $Q$ and $A/P$ are acyclic. \endpf

The lemma is applied above by taking $A= \sing (X(\CC))$ and
    \[ F_pA = \lim_{\stackrel{Y \hookrightarrow X, \text{ closed}}{ \dim Y\leq p}}
    \sing(Y(\CC)) \]
and then $P/Q = \forg \CLD(X)$.

\subsection{Reduction to singular chains}\label{sec:propsofF}

    Since we will work with the d.g.l. $\pp_F$ we tensor every
    Abelian group with $\QQ$ in this section. Thus, for example,
    $\CLD(\bullet)$ really stands for $\CLD(\bullet) \tensor \QQ$, and
    $\sing(\bullet)$ really stands for $\sing(\bullet) \tensor \QQ$.

 We have a diagram
    \begin{equation}
        \xymatrix{
        \fin \times \fin \ar@<-1ex>[d]_{X} \ar@<1ex>[d]^{X} \ar[r]^\amalg & \fin \ar[d]^X \\
        \Aff_k \times \Aff_k \ar@<-1ex>[d]_{\CLD} \ar@<1ex>[d]^{\CLD} \ar[r]^\times & \Aff_k \ar[d]^\CLD \\
        \mm \times \mm \ar[r]^\tensor & \mm
        }
    \end{equation}
    where the top square commutes and the bottom commutes up to natural
    transformation. That is, there is a natural transformation
    \begin{equation}
        N: \CLD(X) \tensor \CLD(Y) \to \CLD(X \times Y)
    \end{equation}
    such that $F = \CLD \circ X^\bullet$ is a functor
    satisfying the axioms of Section \ref{sec:axioms}. We let $\aa
    =\mm$. Since $\mm$ is an Abelian category, all the axioms on $\aa$
    required to apply Theorem \ref{thm:maindgl} are satisfied and the d.g.l.s $\pp_F$ and $\pp_F^N$ make sense.

Nori's construction of $N: \CLD(X) \tensor \CLD(Y) \to \CLD(X \times
Y)$ gives an isomorphism on homology $ H_*(X) \tensor H_*(Y) \cong
H_*(X \times Y)$. This motivates our fifth axiom (continuing the list
from Section \ref{sec:axioms}):
    \begin{enumerate}
        \item[5.] For all $S,T \in \fin$,
                \begin{equation}
                    FS \otimes FT \xrightarrow{N(S,T)} F( S
                    \amalg T)
                \end{equation}
                is a quasi-isomorphism
    \end{enumerate}

Since we want to prove that applying the forgetful functor to $\pp_F$
computes the rational homotopy groups, we wish to show that for the
purposes of computing $H_*(\forg \pp_F)$ or $H_*(\forg \pp_F^N)$,  we
may replace $\forg \CLD$ with $\sing_*$.

By the discussion of the previous section, we have a commutative
diagram of quasi-isomorphisms for any $m,n$
\begin{equation}\label{egn:sig}
    \xymatrix{
        \forg \CLD(X^m)\tensor \forg \CLD(X^n) \ar[d] & P_*(X^m) \tensor P_*(X^n) \ar[l] \ar[r] \ar[d] & \sing_*(X^m) \tensor \sing_*(X^n) \ar[d]  \\
        \forg \CLD(X^{m+n}) & P_*(X^{m+n}) \ar[l] \ar[r] & \sing_*(X^{m+n})
    }
\end{equation}
and this diagram is compatible with the action of $\Sigma_m \times
\Sigma_n$.

Now in this context $\forg \pp^N_F$ is a double complex since each
term $\pp_F(n)$ is a complex itself. If we let $F' = P_*(X^\bullet)$
and $F'' = \sing_*(X^\bullet)$ then the above diagram insures that we
have morphisms of double complexes
\begin{equation}
    \forg \pp_F^N \leftarrow \pp_{F'}^N \to \pp_{F''}^N
\end{equation}
and we have the spectral sequence of a bounded double complex. Taking
the homology from the internal differential (the one other than the
``geometric'' differential) the presence of quasi-isomorphisms in
\eqref{egn:sig} yields an isomorphism at the $E_1$ term and thus
isomorphisms at the $E_\infty$ term. Thus $H_*(\forg \pp_F^N) \cong
H_*(\pp_{F''}^N)$.

Hence we define $F= \sing \circ X^\bullet$ for the rest of this
paper. With this definition the diagram \eqref{eq:star} gives the
geometric cobar construction. By a theorem of Adams \cite{adams},
this complex computes the homology of the loop space of $X$ when $X$
is $1$-connected.

Now $\pp_F$ is a subcomplex of the geometric cobar complex $R$ and
$\pp_F^N$ is a (split) subcomplex of $R_F^N$. Writing this out we
have the commutative diagram:
\begin{equation*}
    \xymatrix{
    \vdots & \vdots \\
    \pp_F(3) = R_{w_3} \sing(X^3) \ar@{^{(}->}[r] \ar[u]^{F(f_3)}  & \sing(X^3)
    \ar[u]^{F(f_3)} \\
    \pp_F(2) = R_{w_2} \sing(X^2) \ar@{^{(}->}[r] \ar[u]^{F(f_2)}  & \sing(X^2)
    \ar[u]^{F(f_2)} \\
    \pp_F(1) = \sing(X) \ar@{^{(}->}[r] \ar[u]^{F(f_1)}  & \sing(X)
    \ar[u]^{F(f_1)} \\
    }
\end{equation*}

In the case of $\pp_F^N$ the top of this diagram would stop at
    \begin{equation*}
        \xymatrix{    \pp_F^N= R_{w_N}\sing(X^N) \ar@{^{(}->}[r]  &
        \sing(X^N)}
    \end{equation*}

\section{Comparison with Rational Homotopy Theory}

\subsection{Basic Definitions}\label{sec:RHBD}

First we recall some facts about the Sullivan complex of polynomial
differential forms with coefficients in $\QQ$. This description is
taken from \cite{fht}.

First we recall the functor from
    \[ \text{Simplicial Sets} \times \text{Simplicial Cochain
    Algebras} \to \text{Cochain Algebras} \]
which assigns to a simplicial set $K$ and a simplicial cochain
algebra $A = \{ A_n \}_{n\geq 0}$ the cochain complex $A(K) = \{
A^p(K)\}_{p\geq 0}$ in the following way. (By a cochain algebra, we
mean a d.g.a. of the form $A = \{ A^p \}_{p\geq 0}$.) $A^p(K)$ is
defined to be the set of simplicial morphisms from $K$ to $A^p$ where
we define, for $\sigma \in K_n$, $\Phi, \Psi \in A^p(K)$, addition by
$(\Phi +\Psi)_\sigma = \Phi_\sigma + \Psi_\sigma$, scalar
multiplication by $(\lambda \Psi)_\sigma = \lambda \Psi_\sigma$ and
the differential by $(d \Psi)_n = d(\Psi_n)$. The algebra structure
is ``point-wise'' multiplication: $(\Psi \cdot \Phi)_\sigma =
\Psi_\sigma \Phi_\sigma$.

Define a simplicial commutative cochain algebra $A_{PL} = \{
(A_{PL})_n \}_{n \geq 0}$ by
    \[ (A_{PL})_n = \frac{\Lambda(t_0, \ldots , t_n, y_0, \ldots ,
    y_n)}{(\sum t_i-1, \sum y_i )} \]
    where $\Lambda$ denotes the free graded commutative algebra (in this case over $\QQ$), the
    $t_i$ are in degree $0$, the $y_i$ are in degree $1$, and the
    differential is defined by $dt_i=y_i$ and $dy_i=0$. The face and
    degeneracy maps are defined by considering $(A_{PL})_n$ as a
    sub-cochain algebra of $A_{DR}(\Delta^n)$, the de-Rham  complex
    on $\Delta^n$, and using the face and degeneracy maps induced from the
    simplices $\Delta^n$.

Now from a topological space $X$ we use the standard way of obtaining
a simplicial set. $S_n(X)$ is defined to be the set of all singular
$n$-simplices $\sigma: \Delta^n \to X$, where the $i$-th face map,
$\partial_i: S_{n+1}(X)  \to S_n(X)$ is given by pulling back along
the face inclusions $\Delta^n \hookrightarrow \Delta^{n+1}$ and the
$i$-th degeneracy is given by the map from $\Delta^{n+1}$ to
$\Delta^n$ that collapses the $j$-th face.  Then applying the above
construction with the simplicial set $\{S_n(X)\}$ and the simplicial
cochain algebra $A_{PL}$ gives the  commutative cochain algebra of
polynomial rational forms on $X$, $A_{PL}(X)$. This is a functor from
spaces to commutative cochain algebras.

If we choose a point $p \in X$ the inclusion $p \hookrightarrow X$
determines map $A_{PL}(X) \to A_{PL}(p)= \QQ$ which makes $A_{PL}(X)$
an augmented algebra.

Let us recall some facts about the bar construction. Let  $(A, m,
\epsilon)$ be an augmented d.g.a. over $k$, $A \xrightarrow{\epsilon}
k$. Let $IA= \ker(\epsilon)$, and left $s$ denote the shift
(suspension) map. Then $\BAR(A)$ is a coaugmented d.g.c. As a
coaugmented graded coalgebra it is the tensor
    coalgebra
        \[ T(s (IA)). \]
        That is, the  coalgebra structure is given by
        \[ [a_1 | \cdots | a_n] \mapsto [a_1 | \cdots | a_n] \otimes
        1 + \sum_{i=1}^{n-1}[a_1| \cdots | a_i] \tensor [a_{i+1}|
        \cdots | a_n] + 1 \tensor [a_1 | \cdots | a_n]. \]
    The differential (which we will not write out here) can be
    expressed
        \[ d_B = d_I + d_E \]
    where the internal differential, $d_I$, comes from the
    differential on $A$, while the external differential $d_E$ comes
    from the multiplicative structure in $A$.  $\BAR(A)$ may be
    considered a double complex by indexing with $-s$ the tensor
    (external) degree and with $t$ the differential form (internal)
    degree,
    \begin{equation}
        B^{-s,t}(A) = [ \otimes^sIA]^t.
    \end{equation}
    The negative degree allows us to use the total degree $t-s$ as
    the cohomolgical degree of the total complex $\BAR(A)$.

    If we {\it forget the differential} then $T(s(IA))$ has a
    product, the shuffle product,
        \[ sh: T(s(IA)) \tensor T(s(IA)) \to T(s(IA)).  \]
    For any d.g. vector space, the shuffle
    product exists and is compatible with the tensor coalgebra
    structure.

    In general the shuffle product is not compatible with
    the differential.  It is compatible when $(A,m,\epsilon)$ is
    graded commutative, in which case $\BAR(A)$ become a d.g. Hopf
    algebra. We will always use commutative algebras.

    Forgetting the external differential, for any d.g. vector
    space  $V$ we can define $\BART(V)$ as the tensor
    coalgebra $T(s(V))$ with the internal differential and the
    shuffle product. This is also a d.g. Hopf algebra. Notice that
    except for the differential, $\BART(IA)$ and $\BAR(A)$ are the
    same, in particular they have the same tensor coalgebra
    structure.
    (This is equivalent to $\BAR(A)$ when $A$ has all products equal
    to zero.) 

Furthermore $\BAR(A)$ has the ``external'' filtration (called the
``bar'' filtration) defined by $\mathfrak{B}^{-s} = \oplus_{u \leq s}
B^{-u,v}$. So $\QQ= \mathfrak{B}^0 \subset \mathfrak{B}^{-1} \subset
\ldots \subset \BAR(A)$. Then projection onto $\mathfrak{B}^0$
defines an augmentation which make $\BAR(A)$ into an augmented d.g.
Hopf algebra.

Since $\BAR(A)$ or $\BART(A)$ is an augmented d.g. Hopf algebra, we
can form the Lie coalgebra of indecomposables: Let $J$ be the
augmentation ideal of $\BAR(A)$. Then the cokernel of (shuffle)
multiplication on $J$ defines the indecomposables:
    \begin{equation}
        J \tensor J \xrightarrow{sh} J \to J/J^2 = Q\BAR(A) \to 0
    \end{equation}

Recall that a space $X$ is nilpotent if for all $p\in X$,
$\pi_1(X,p)$ is a nilpotent group and the natural action of
$\pi_1(X,p)$ on $\pi_i(X,p) \tensor \QQ$ is unipotent for $i \geq 2$.

In \cite{hain1} Hain defines the homotopy Lie algebra
$\mathfrak{g}_*(X,p)$ to be the graded Lie algebra where
$\mathfrak{g}_0$ is the $\QQ$-Malcev completion of $\pi_1(X,p)$ and
for $i\geq 2$ define $\mathfrak{g}_i(X,p) = \pi_{i+1}(X,p) \tensor
\QQ$ if $X$ is nilpotent, and $0$ otherwise. This definition is not
necessarily compatible with the modern topology literature, so we
will call this the ``homotopy Lie algebra in the sense of Hain.''

Now let $A(X)=A_{PL}(X)$ be the $\QQ$--de Rham complex as in
\cite{hain1}. Hain proves that if $E^* \to \QQ$ is an augmented
commutative d.g. algebra which is quasi-isomorphic to $A(X)$ with
augmentation given by a point $p$, then if $X$ is $1$-connected,
integration  induces a natural Lie coalgebra isomorphism of the dual
of the homotopy Lie algebra (in the sense of Hain) and the homology
of the indecomposables of the bar construction on $E^*$:
    \begin{equation}
        (\mathfrak{g}_*(X,p))^\vee \cong H_*( Q\BAR(E^*))
    \end{equation}
If $X$ is not simply connected, then we only assert that
    \begin{equation}
        (\mathfrak{g}_0(X,p))^\vee \cong H_0( Q\BAR(E^*)).
    \end{equation}

 This will serve as the foundation of our comparison theorem.

Notice with $X = (X,p)$, working with $\sing(X)$,
    \begin{equation}
    0 \to \sing(p) \to \sing(X) \to \sing(X,p) \to 0
    \end{equation}
 is dual to working
with $IA$
    \begin{equation}
    0 \to IA \to A(X) \to A(p) \to 0.
    \end{equation}
Furthermore we will tensor with the rationals throughout, that is
$\sing = \sing \tensor \QQ$.

\subsection{Alternate description of the bar construction}

We can give another description of the bar construction from the
point of view of sections 2 and 3. That is, rather than considering
functors  on the category $\dd$, we recast the bar construction in
terms of the functor corepresented by a commutative algebra.

For this section, let $\aa$ be an additive $\tensor$ category
satisfying the usual commutativity, associativity and unity
conditions. Fix an object $A \in \aa$ and assume that $m: A \tensor A
\to A$ is a map satisfying the constraints of commutative
multiplication.

Then there is a functor $A^{\tensor \bullet}$ from $\dd$ to $\aa$
given by
    \[ [n] \mapsto A^{\tensor n} \]
    and on morphisms, if two elements map to the same element, the
    functor multiplies those terms. This is well defined since $A$ is
    a commutative algebra. For example the map $[2] \to [1]$ induces
    multiplication $A^{\tensor 2} \xrightarrow{m} A$. On the other
    hand if a map in $\dd$ misses an element in the target, the
    unit object of $A$ is inserted. For example the map $[1] \to [2]$
    which sends the singleton of $[1]$ to the first element of $[2]$
    induces the map $A \to A^{\tensor 2}$,  $a \mapsto a \tensor 1$.

Then the algebraic bar construction on $A$ is a functor
$\BAR^A(\bullet) : \dd \to \Ch(\aa)$ obtained by applying $A^{\tensor
\bullet}$ to \eqref{simpcobar}. Of course $\BAR^A_{-N}$ is given by
applying $A^{\tensor \bullet}$ to \eqref{simpcobarN}. Re-interpreting
the work of Section 3 in this context we immediately get that
$\BAR^A$ is a differential graded coalgebra, and thus has an
underlying Lie coalgebra structure. Henceforth define $A = A(X)$, the
$\QQ$--de Rham complex,  as in Section 5.1, for a fixed $X$. The plan
of our comparison is to use the natural transformation of functors
from $\fin$ to $\Ch(\Ab)$ given by integration
    \[ \int : (IA)^{\tensor \bullet} \to \sing(X^\bullet)^\vee. \]
Integration here makes sense as we are have differential forms which
can be integrated over simplices.

\subsection{Comparison}

Since the shuffle product is a map of degree $0$, we have that the
cokernel of $\mathfrak{B}^{-N} \tensor \mathfrak{B}^{-N}
\xrightarrow{sh} \mathfrak{B}^{-N}$, which might be called
$Q\mathfrak{B}^{-N}$ is the same as taking all terms of $Q\BAR(A)$ of
external degree at smallest $-N$, that is \[Q\mathfrak{B}^{-N} =
\bigoplus_{-s>-N}Q\BAR^{-s}(A). \]

Since for each $n \in \NN$ integration gives a map
    \begin{equation}
        IA^{\tensor n} \xrightarrow{\int} \sing(X^n)^\vee
    \end{equation}
we will show this will give a map
    \begin{equation}
        Q\mathfrak{B}^{-N} \to (\pp^N_F \tensor \QQ)^\vee.
    \end{equation}

\begin{thm}\label{thm:comparison} Let $F$ be the functor given by
$\forg \CLD(X^\bullet) \tensor \QQ$, then $(\pp_F^N)^\vee$ is a Lie
coalgebra such that
    \[   H^*( Q \mathfrak{B}^{-N} ) \cong H^*( (\pp_F^N)^\vee)   \] as Lie coalgebras.
\end{thm}

\pf Since $F$  satisfies
    the conditions of Section \ref{sec:propsofF}, the statement of the theorem
    reduces to the case $F= \sing(X^\bullet) \tensor \QQ$. Notice that the construction
    of $\pp_F$  does not use
    the Alexander--Whitney comultiplication. Since the underlying
    coalgebra structure on the bar construction comes from the free
    tensor coalgebra, and the coalgebra structure on
    $(\pp_F^N)^\vee$ arises from the dual of the free tensor algebra
    structure, their coalgebra structures are automatically
    compatible.
In this proof, after unravelling the definition of $F$, we know that
the image of $F(\widetilde{w_n})$, $F(n) F(\widetilde{w_n})$,  is
identified with $R_{w_n}\sing(X^n)$. Let us denote the dual of the
map $F(\widetilde{w_n}): \sing(X^n) \to \sing(X^n)$ by $w_n^\vee :
\sing(X^n)^\vee \to \sing(X^n)^\vee$ so that we may write the image
of this second map as $w_n^\vee \sing(X^n)^\vee$.

We need several lemmas to complete this proof:

\begin{lemma}
    \[ w_n^\vee A^{\tensor n} = (Q \BAR (A))^n. \]
\end{lemma}
\pf
    Since the assertion is independent of the differential, we have
    to show that for any vector space $V$,
        \[ w_n^\vee( V^{\otimes n}) = Q \BART(V)^n  \]

    First let $W$ be a finite dimensional vector space.
    By Lemma \ref{lemmaquillen}, the following sequence is exact (where $w_n$ splits $i$).
    \begin{equation}
        \xymatrix{
        0 \ar[r] & L_n(W^\vee) \ar[rr]_{i} & & \ar@/_/[ll]_{w_n} (W^\vee)^{\tensor n} \ar[rr]^{n-w_n} & & (W^\vee)^{\tensor n}
        }
    \end{equation}
    where $L_n$ refers to the degree $n$ part of the graded free Lie
    algebra on $W^\vee$ concentrated in degree one.
    Applying $\Hom(- , k)$, we have the exact sequence
    \begin{equation}
        \xymatrix{
        (W)^{\tensor n} \ar[rr]^{(n-w_n)^\vee} & & (W)^{\tensor n}  \ar[rr]_{i^\vee}& &
        \ar@/_/[ll]_{w_n^\vee} (L_n(W^\vee))^\vee \ar[r] & 0.
        }
    \end{equation}
    Dualizing the degree $n$ term of the free Lie algebra on $W^\vee$
    gives the degree $n$ term of the free Lie coalgebra on $W$.

    Meanwhile, $Q(\BART(W^\vee))^n$ is the dual of $\prim(W)_n$
    inside the tensor algebra $T(W)$. Considering the discussion in
    Section 2.4, identifying the free Lie algebra on $W$ with the
    primitive Lie algebra inside the tensor algebra on $W$,    we
    know that  $\prim(W)_n = L_n(W)$,  and so the Lemma is proved for
    $W$.

    Now consider any $V$ and let $W_1 \subset W_2 \subset \ldots
    \subset W$ be an exhaustive filtration of finite dimensional subspaces. For each $i$
    we have $w_n^\vee W_i^{\tensor n} = (Q \BAR (W_i))^n$.
    Furthermore, for $i<j$ we have
    \begin{equation}
    \xymatrix{
        w_n^\vee W_i^{\tensor n} \ar@{=}[r] \ar@{^{(}->}[d] & (Q \BART (W_i))^n \ar@{^{(}->}[d] \\
        w_n^\vee W_j^{\tensor n} \ar@{=}[r] & (Q \BART (W_j))^n
    }
    \end{equation}
    since a basis for $W_i$ can be completed to a basis of $W_j$
    and since the relations determined by the shuffle product applied
    to $\BART(W_j)$ do not kill any generators of $(Q \BART(W_i))^n$.
    Since a vector space is the colimit of an exhaustive filtration,
    the Lemma follows. \endpf

\begin{lemma}\label{littlelemma}
    If $m:M_* \to M_*$ and $n:N_* \to N_*$ are morphisms of chain
    complexes which are idempotent, and $f: M_* \to N_*$ is a
    quasi-isomorphism
    which makes
\[            \xymatrix{
            M_* \ar[r]^f \ar[d]^m & N_* \ar[d]^n \\
            M_* \ar[r]^f & N_*
            } \]
    commute, then $\im(m) \xrightarrow{f} \im(n)$ is a
    quasi-isomorphism.
\end{lemma}
\pf
    Since $m$ and $n$ are idempotent, each provides a splitting of complexes,
    $M \cong mM \oplus (1-m)M$ and $N \cong nN \oplus (1-n)N$. Use
    these isomorphisms as the vertical arrows in the commutative
    diagram
    \[ \xymatrix{
        M \ar[r]^f \ar[d] & N \ar[d] \\
        mM \oplus (1-m)M \ar[r]^{f \oplus f} & nN \oplus (1-n)N
        } \]
    Since homology commutes with direct sum, we conclude that
    $H_*(mM) \cong H_*(nN)$. \endpf

Let \[ e_n = \frac{w_n}{n}: V^{\otimes n} \to V^{\otimes n}. \] Then
$e_n$ is an idempotent, $e_n^2=e_n$.

\begin{lemma}\label{lem:intisquasi}
    Integration, $\int: w_n^\vee (IA)^{\tensor n}\to w_n^\vee
    \sing(X^n)^\vee$ is a quasi-isomorphism.
\end{lemma}
    \pf First notice that for $\sigma \in \Sigma_n$, the following
    diagram commutes, and $\int$ is a quasi isomorphism.
    \begin{equation}
        \xymatrix{
        IA^{\tensor n} \ar[r]^{\int} \ar[d]^\sigma & \sing(X^n)^\vee
        \ar[d]^\sigma \\
        IA^{\tensor n} \ar[r]^{\int} & \sing(X^n)^\vee
        }
    \end{equation}
    The result now follows immediately from Lemma \ref{littlelemma} by
    using $e_n$ in place of $w_n$ and observing that for vector spaces,
    $e_n$ and $w_n$ have the same image. \endpf

Now expanding out $\BAR(A)$ and the dual of the geometric cobar
construction along their external degrees, and integrating in each
external degree we have a commutative diagram
    \begin{equation}\label{ninetynine}
        \xymatrix{
            IA \ar[d] & IA \tensor IA \ar[l] \ar[d] & IA^{\tensor 3}
            \ar[d] \ar[l] & \cdots \ar[l] \\
            \sing(X)^\vee & \sing(X \times X)^\vee \ar[l] &
            \sing(X^3)^\vee \ar[l] & \ldots \ar[l]
        }
    \end{equation}
Now we wish to show that applying the $w_{i}^\vee$ to each term gives
us a subcomplex:

\begin{lemma} \label{lemma:setupcomparison}
    The following diagram commutes:
\[    \xymatrix{
        w_{n+1}^\vee\sing(X^{n+1})^\vee \ar[r] &
        w_{n}^\vee\sing(X^{n})^\vee \ar[r] &
        w_{n-1}^\vee\sing(X^{n-1})^\vee \\
        w_{n+1}^\vee IA^{\tensor n+1} \ar[r] \ar[u]^{\int} &
        w_{n}^\vee IA^{\tensor n} \ar[r] \ar[u]^{\int} &
        w_{n-1}^\vee IA^{\tensor n-1} \ar[u]^{\int}
    } \]
\end{lemma}
\pf  
    As pointed out the the end of Section 5.2, since integration is a
    natural transformation of functors
    \[ \int : (IA)^{\tensor \bullet} \to \sing(X^\bullet)^\vee \]
    which respects the $\Sigma_n$ action the lemma follows from
    \eqref{ninetynine}.
\endpf

    Now by Lemma
    \ref{lemma:setupcomparison} we have a morphism of double complexes
    from $Q \mathfrak{B}^{-N}$ to $(\pp_F^N)^\vee$ which is an isomorphism on the
    $E_1$ term of the spectral sequence of the double complex. Since both of these
    complexes are bounded (by $N+1$) in external degree, this
    implies that they have the same $E_\infty$ term, namely
    $H^*(Q\mathfrak{B}^{-N}) \cong H^*((\pp_F^N)^\vee)$. As stated before both double complexes have
    comultiplication from the tensor coalgebra structure, thus this
    isomorphism is an isomorphism of Lie coalgebras. This proves Theorem \ref{thm:comparison}. \endpf

\begin{cor} If $X$ is $1$-connected then
    \[ \lim_{\overleftarrow{N}}H^*((\pp_F^N)^\vee) \cong
    (\mathfrak{g}_*)^\vee, \] and if $X$ is only assumed to be connected then
    \[ \lim_{\overleftarrow{N}}H^0((\pp_F^N)^\vee) \cong
    (\mathfrak{g}_0)^\vee \]
    as Lie coalgebras and the cobracket is the
    dual of the Whitehead bracket.
\end{cor}
\pf This follows from Hain's Theorem 2.6.2 in \cite{hain1}
identifying the homology of the indecomposables of the bar
construction with the dual of rational homotopy Lie algebra.
\endpf

\begin{cor}
    If $X$ is of
    finite type (e.g. if $X$ is an algebraic variety) and if $X$ is $1$-connected, then for fixed $k$ there exists some $N>0$ such that
    \[ H_i( F(\pp^N)^\vee) \xrightarrow{\cong} (\pi_{i+1}(X,p)\tensor \QQ)^\vee \qquad \forall \, i \leq k \]
\end{cor}
    \pf Look at the $E_1$ term of the spectral sequence of the bar
    construction obtained by taking the homology in the internal
    degree. For any $q \geq 2$ the terms that contribute to the
    calculation of $\pi_q(X,p) \tensor \QQ$ sit in total degree
    $q-1$, and thus are $H^q(X) \oplus H^{q+1}(X^2) \oplus \ldots
    \oplus    H^{q+j}(X^j) \oplus \ldots$.  But the Kunneth Theorem
    tells us that
    \[ H^{q+j}(X^j) = \bigoplus_{i_1 + \ldots + i_j=
    q+j}H^{i_1}(X)\tensor \ldots \tensor H^{i_j}(X) \]
    and once $j>q$ then at least one of the $i_k=1$ and the simply
    connected hypothesis forces that term, and all higher terms to be
    zero. \endpf

\end{document}